\documentclass[12pt,english]{elsarticle}
\usepackage[T1]{fontenc}
\usepackage[latin9]{inputenc}
\usepackage{geometry}
\usepackage{babel}
\usepackage{amsmath}
\usepackage{amssymb}
\usepackage{graphicx}
\usepackage{esint}
\usepackage[unicode=true]
 {hyperref}

\makeatletter

\providecommand{\tabularnewline}{\\}

\usepackage[font=footnotesize,labelfont=bf]{caption}
\usepackage{lineno}

\usepackage{sectsty}
\sectionfont{\large}
\subsectionfont{\normalsize\it\bfseries}
\subsubsectionfont{\normalsize\it}

\makeatother

\begin{document}

\title{\textbf{\large{A consistent operator splitting algorithm and a two-metric
variant: Application to topology optimization}}}

\author{Cameron Talischi$^{*}$, Glaucio H. Paulino\medskip{}
}

\address{Department of Civil and Environmental Engineering, University of
Illinois at Urbana-Champaign, USA\smallskip{}
}

\address{$^{*}$\emph{Corresponding author, \href{mailto:ktalisch@illinois.edu}{ktalisch@illinois.edu}}}
\begin{abstract}
In this work, we explore the use of operator splitting algorithms
for solving regularized structural topology optimization problems.
The context is the classical structural design problems (e.g., compliance
minimization and compliant mechanism design), parameterized by means
of density functions, whose ill-posendess is addressed by introducing
a Tikhonov regularization term. The proposed forward-backward splitting
algorithm treats the constituent terms of the cost functional separately
which allows suitable approximations of the structural objective.
We will show that one such approximation, inspired by the optimality
criteria algorithm and reciprocal expansions, improves the convergence
characteristics and leads to an update scheme that resembles the well-known
heuristic sensitivity filtering method. We also discuss a two-metric
variant of the splitting algorithm that removes the computational
overhead associated with bound constraints on the density field without
compromising convergence and quality of optimal solutions. We present
several numerical results and investigate the influence of various
algorithmic parameters. 

\smallskip{}

\end{abstract}
\maketitle

\section{Introduction}

\thispagestyle{plain}The goal of topology optimization is to find
the most efficient shape of a physical system whose behavior is captured
by the solution to a boundary value problem that in turn depends on
the given shape. As such, optimal shape problems can be viewed as
a class of optimal control problems in which the control is the shape
or domain of the governing state equation. These problems are in general
ill-posed in that they do not admit solutions in the classical sense.
For example, the basic compliance minimization problem in structural
design, wherein one aims to find the stiffest arrangement of a fixed
volume of material, favors non-convergent sequences of shapes that
exhibit progressively finer features (see, for example, \citep{Allaire-book}
and reference therein). A manifestation of the ill-posedness of the
continuum problem is that naive finite element approximations of the
problem may suffer from numerical instabilities such as spurious checkerboard
patterns or exhibit mesh-dependency of the solutions, both of which
can be traced back to the absence of an internal length-scale in the
continuum description of the problem \citep{Sigmund:1998p441}. An
appropriate regularization scheme, based on one's choice of parametrization
of the unknown geometry, must therefore be employed to exclude this
behavior and limit the complexity of the admissible shapes. 

One such restriction approach, known as the \emph{density filtering}
method, implicitly enforces a prescribed degree of smoothness on all
the admissible density fields that define the topology \citep{Bourdin:2001p457,Bruns:2001p2341}.
This method and its variations are consistent in their use of sensitivity
information in the optimization algorithm since the sensitivity of
the objective and constraint functions are computed with respect to
the associated auxiliary fields whose filtering defines the densities%
\footnote{Effectively filtering is a means to describe the space of admissible
densities with an embedded level of regularity -- for more refer to
\citep{PolyTop}.%
}. By contrast, the \emph{sensitivity filtering} method \citep{Sigmund:1998p441,SigmundMaute2013},
which precedes the density filters and is typically described at the
discrete level, performs the smoothening operation directly on the
sensitivity field after a heuristic scaling step. The filtered sensitivities
then enter the update scheme that evolves the design despite the fact
they do not correspond to the cost function of the optimization problem.
While the sensitivity filtering has proven effective in practice for
certain class of problems (for compliance minimization, it enjoys
faster convergence than the density filter counterpart), a proper
justification has remained elusive. As pointed out by Sigmund \citep{Sigmund:2007p460},
it is generally believed that ``the filtered sensitivities correspond
to the sensitivities of a smoothed version of the original objective
function'' even though ``it is probably impossible to figure out
what objective function is actually being minimized.'' This view
is confirmed in the present work, as we will show that an algorithm
with calculations similar to what is done in the sensitivity filtering
can be derived in a \emph{consistent} manner from a proper regularization
of the objective.

The starting point is the authors' recent work \citep{Talischi-JMPS}
on an operator splitting algorithm for solving the compliance minimization
problem where a Tikhonov regularizaton term is introduced to address
the inherent ill-posedness of the problem. The derived update expression
naturally contains a particular use of Helmholtz filtering, where
in contrast to density and sensitivity filtering methods, the filtered
quantity is the gradient descent step associated with the original
structural objective. The key observation made here is that if the
gradient descent step in this algorithm is replaced by the optimality
criteria (OC) update, then the interim density has a similar form
to that of the sensitivity filter and in fact produces similar results
(cf. Figure \ref{fig:ResultSenOCvsModFBS}). To make such a leap rigorous,
we essentially embed the same reciprocal approximation of compliance
that is at the heart of the OC scheme in the forward-backward algorithm.
This leads to a variation of the forward-backward splitting algorithm
in \citep{Talischi-JMPS} that is consistent, demonstrably convergent
and computationally tractable.

Within the more general framework presented here, we will examine
the choice of move limits and the step size parameter more closely
and discuss strategies that can improve the convergence of the algorithm
while maintaining the quality of final solutions. We also discuss
a two-metric variant of the splitting algorithm that removes the computational
overhead associated with the bound constraints on the density field
without compromising convergence and quality of optimal solutions.
In particular, we present and investigate scheme based on the two-metric
projection method of \citep{Bert-ProjNewton,Gafni-TMP} that allows
for the use of a more convenient metric for the projection step enforcing
these bound constraints. This algorithm requires a simple and computationally
inexpensive modification to the splitting scheme but features a min/max-type
projection operation similar to OC-based filtering methods. We will
see from the numerical examples that the two-metric variation retains
the convergence characteristics of the forward-backward algorithm
for various choices of algorithmic parameters. The details of the
two types of algorithms are described for the finite-dimensional optimization
problem obtained from the usual finite element approximation procedure,
which we prove is convergent for Tikhonov-regularized compliance minimization
problem.

The remainder of this paper is organized as follows. In the next section,
we describe the model topology optimization problem and its regularization.
A general iterative scheme---one that encompasses the previous work
\citep{Talischi-JMPS}---for solving this problem based on forward-backward
splitting scheme is discussed in section 3. Next, in section 4, the
connection is made with the sensitivity filtering method and the OC
algorithm, and the appropriate choice of the approximate Hessian is
identified. For the sake of concision and clarity, the discussion
in these three sections is presented in the continuum setting. In
section 5, we begin by showing that the usual finite element approximations
of the Tikhonov-regularized compliance minimization problem are convergent
and derive the vector form of the discrete problem. The proposed algorithms
along with some numerical investigation are presented in sections
6 and 7. We conclude the work with some closing remarks and future
research directions in the section 8.

Before concluding the introduction, we briefly describe the notation
adopted in this paper. As usual, $L^{p}(\Omega)$ and $H^{k}(\Omega)$
denote the standard Lebesgue and Sobolev spaces defined over domain
$\Omega$ with their vector-valued counterparts $L^{p}(\Omega;\mathbb{R}^{d})$
and $H^{k}(\Omega;\mathbb{R}^{d})$, and $L^{p}(\Omega;K)=\left\{ f\in L^{p}(\Omega):f(\mathbf{x})\in K\mbox{ a.e.}\right\} $
for a given $K\subseteq\mathbb{R}$. Symbols $\wedge$ and $\vee$
denote the point-wise min/max operators. Of particular interest are
the inner product and norm associated with $L^{2}(\Omega)$, which
are written as $\left\langle \cdot,\cdot\right\rangle $ and $\left\Vert \cdot\right\Vert $,
respectively. Similarly, the inner product, norm and semi-norm associated
with $H^{k}(\Omega)$ are denoted by $\left\langle \cdot,\cdot\right\rangle _{k}$,
$\left\Vert \cdot\right\Vert _{k}$ and $\left|\cdot\right|_{k}$,
respectively. Given a bounded and positive-definite linear operator
$\mathcal{\mathcal{B}}$, we write $\left\langle u,v\right\rangle _{\mathcal{\mathcal{B}}}\equiv\left\langle u,\mathcal{\mathcal{B}}v\right\rangle $
and the associated norm by $\left\Vert u\right\Vert _{\mathcal{B}}\equiv\left\langle u,u\right\rangle _{\mathcal{B}}^{1/2}$.
Similarly, the standard Euclidean norm of a vector $\mathbf{v}\in\mathbb{R}^{m}$
is denoted by $\left\Vert \mathbf{v}\right\Vert $ and given a positive-definite
matrix $\mathbf{B}$, we define $\left\Vert \mathbf{v}\right\Vert _{\mathbf{B}}=\left(\mathbf{v}^{T}\mathbf{B}\mathbf{v}\right)^{1/2}$.
The $i$th components of vector $\mathbf{v}$ and the $(i,j)$-th
entry of matrix $\mathbf{B}$ are written as $\left[\mathbf{v}\right]_{i}$
and the $\left[\mathbf{B}\right]_{ij}$, respectively.

\section{Model Problem and Regularization}

\begin{figure}
\begin{centering}
\includegraphics[scale=0.55]{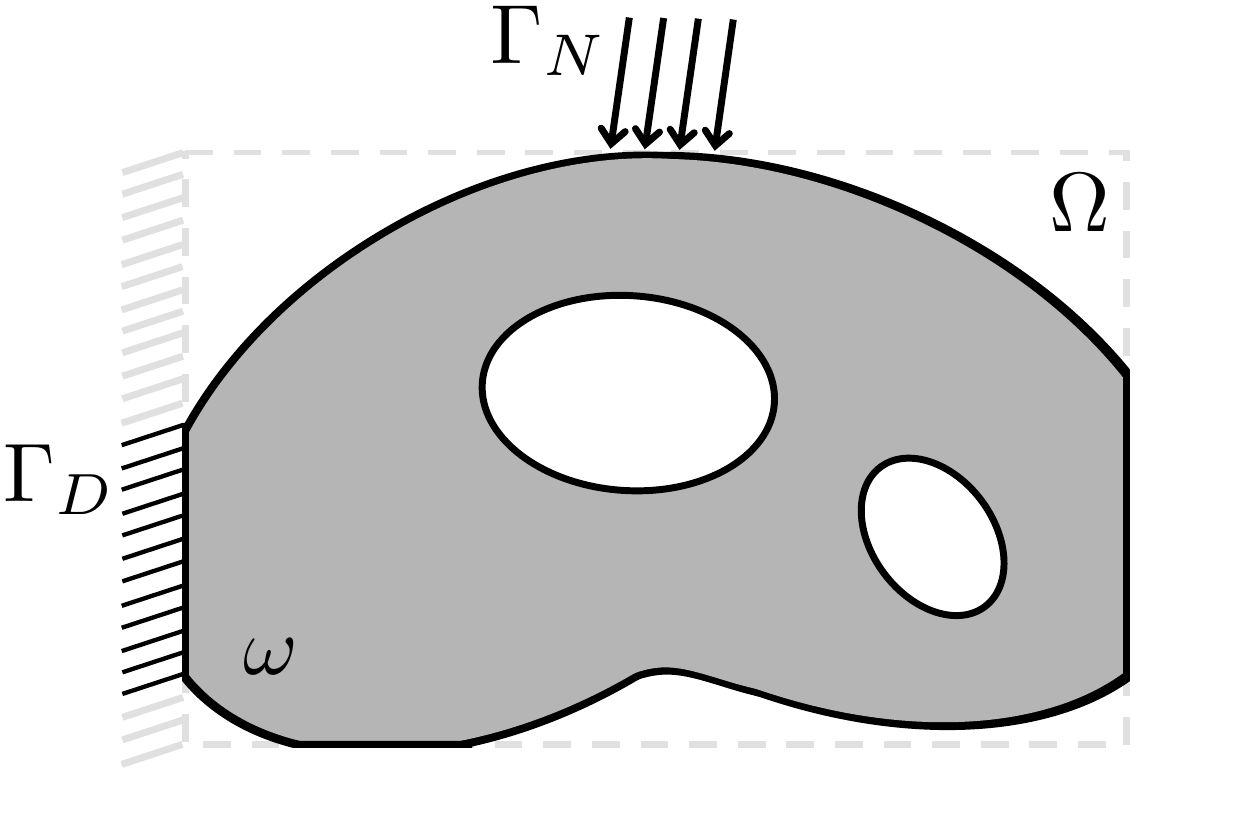}
\par\end{centering}

\caption{Illustration of the prescribed boundary conditions defined on the
design domain $\Omega$. In a density formulation, each admissible
shape $\omega\subseteq\Omega$ can be associated with some density
function $\rho\in L^{\infty}(\Omega;\left[\delta_{\rho},1\right])$\label{fig:DomainSchematics}}
\end{figure}

We begin with the description of the compliance minimization problem
which is used as the model problem in this work. Let $\Omega\subseteq\mathbb{R}^{d},d=2,3$
be the extended design domain with sufficiently smooth boundary. We
consider boundary segments $\Gamma_{D}$ and $\Gamma_{N}$ that form
a nontrivial partition of $\partial\Omega$, i.e., $\Gamma_{D}\cap\Gamma_{N}=\emptyset$,
$\partial\Omega=\overline{\Gamma}_{D}\cup\overline{\Gamma}_{N}$ and
$\Gamma_{D}$ has non-zero surface measure (see Figure \ref{fig:DomainSchematics}).
Each design over $\Omega$ is represented by a density function $\rho$
whose response is characterized by the solution $\mathbf{u}_{\rho}$
to the elasticity boundary value problem, given in the weak form by
\begin{equation}
a(\mathbf{u},\mathbf{v};\rho)=\ell(\mathbf{v}),\qquad\forall\mathbf{v}\in\mathcal{V}\label{eq:weak-form}
\end{equation}
where $\mathcal{V}=\{\mathbf{u}\in H^{1}(\Omega;\mathbb{R}^{d}):\mathbf{u}|_{\Gamma_{D}}=\mathbf{0}\}$
is the space of admissible displacements and
\begin{equation}
a(\mathbf{u},\mathbf{v};\rho)=\int_{\Omega}\rho^{p}\mathbf{C}_{0}\boldsymbol{\epsilon}(\mathbf{u}):\boldsymbol{\epsilon}(\mathbf{v})\mathrm{d}\mathbf{x},\qquad\ell(\mathbf{v})=\int_{\Gamma_{N}}\mathbf{t}\cdot\mathbf{v}\mathrm{d}s\label{eq:a-ell-defs}
\end{equation}
are the usual energy bilinear and load linear forms. Moreover, $\boldsymbol{\epsilon}(\mathbf{u})=(\nabla\mathbf{u}+\nabla\mathbf{u}^{T})/2$
is the linearized strain tensor, $\mathbf{t}\in L^{2}(\Gamma_{N};\mathbb{R}^{d})$
is the prescribed tractions on $\Gamma_{N}$ and $\mathbf{C}_{0}$
is the elasticity tensor for the constituent material. Observe that
the classical Solid Isotropic Material with Penalization (SIMP) model
is used to describe the dependence of the state equation on the density
field, namely that the stiffness is related to the density through
the power law relation $\rho^{p}$ \citep{Bendsoe:1989p447,Rozvany:1992p981,Rozvany:2009p2274}%
\footnote{We use the classical SIMP parametrization with a positive lower bound
on the densities. The reason is that later, we will consider Taylor
expansions in $1/\rho$.%
}. The bilinear form is continuous and also coercive provided that
$\rho$ is measurable and bounded below by some small positive constant
$0<\delta_{\rho}\ll1$. In fact, there exist positive constants $c$
and $M$ such that for all $\rho\in L^{\infty}(\Omega;\left[\delta_{\rho},1\right])$,
\begin{equation}
\left|a(\mathbf{u},\mathbf{v};\rho)\right|\leq M\left\Vert \mathbf{u}\right\Vert _{1}\left\Vert \mathbf{v}\right\Vert _{1},\quad a(\mathbf{u},\mathbf{u};\rho)\geq c\left\Vert \mathbf{u}\right\Vert _{1}^{2},\qquad\forall\mathbf{u},\mathbf{v}\in\mathcal{V}\label{eq:coer-cont-of-a}
\end{equation}
Together with continuity of the linear form $\ell$ (which follows
from the assumed regularity of the applied tractions), these imply
that (\ref{eq:weak-form}) admits a unique solution $\mathbf{u}_{\rho}$
for all $\rho\in L^{\infty}(\Omega;\left[\delta_{\rho},1\right])$.
Moreover, we have the uniform estimate $\left\Vert \mathbf{u}_{\rho}\right\Vert _{1}\leq c^{-1}\left\Vert \mathbf{t}\right\Vert $.
For future use, we also recall that by the principle of minimum potential,
$\mathbf{u}_{\rho}$ is characterized by 
\begin{equation}
\mathbf{u}_{\rho}=\underset{\mathbf{v}\in\mathcal{V}}{\mbox{argmin}}\left[\frac{1}{2}a(\mathbf{v},\mathbf{v};\rho)-\ell(\mathbf{v})\right]\label{eq:min-pot_energy}
\end{equation}
where the term in the bracket is the potential energy associated with
deformation field $\mathbf{v}$. The following is a result that will
be used later in the paper and readily follows from the stated assumptions
(see, for example, \citep{Borrvall:2001p422}): Given a sequence $\left\{ \rho_{n}\right\} $
and $\rho$ in $L^{\infty}(\Omega;\left[\delta_{\rho},1\right])$
such that $\rho_{n}\rightarrow\rho$ strongly in $L^{p}(\Omega),1\leq p\leq\infty$,
the associate displacement fields $\mathbf{u}_{\rho_{n}}$, up to
a subsequence, converge in the strong topology of $H^{1}(\Omega;\mathbb{R}^{d})$
to $\mathbf{u}_{\rho}$. This shows that if the cost functional depends
continuously on $\left(\rho,\mathbf{u}\right)$ in the strong topology
of $L^{p}(\Omega)\times H^{1}(\Omega;\mathbb{R}^{d})$, then compactness
of the space of admissible densities in $L^{p}(\Omega)$ is a sufficient
condition for existence of solutions.

The cost functional for the compliance minimization problem is given
by
\begin{equation}
J(\rho)=\ell(\mathbf{u}_{\rho})+\lambda\int_{\Omega}\rho\mathrm{d}\mathbf{x}\label{eq:compliance}
\end{equation}
The first term in $J$ is the compliance of the design while the second
term represents a penalty on the volume of the material used. Minimizing
this cost functional amounts to finding the stiffest arrangement while
using the least amount of material with elasticity tensor $\mathbf{C}_{0}$.
The parameter $\lambda>0$ determines the trade-off between the stiffness
provided by the material and the amount that is used (which presumably
is proportional to the cost of the design). Since the SIMP model assigned
smaller stiffness to the intermediate densities compared to the their
contribution to the volume, it is expected that in the optimal regime,
the density function are nearly binary (taking only values of $\delta_{\rho}$
and 1) provided that the penalty exponent $p$ is sufficiently large.

As discussed in the introduction, the compliance minimization problem
does not admit a solution unless additional restrictions are placed
on the regularity of density functions. This may be accomplished by
addition of a Tikhonov regularization term to the cost function \citep{Borrvall:2001p2519,Talischi-JMPS}:
\begin{equation}
\min_{\rho\in\mathcal{A}}\ \tilde{J}(\rho)=J(\rho)+\frac{\beta}{2}\left|\rho\right|_{1}^{2}\label{eq:reg-compliance}
\end{equation}
where $\beta>0$ is a positive constant determining the influence
of this regularization (larger $\beta$ leads to smoother densities
in the optimal regime). The minimization of $\tilde{J}$ is carried
out over the set of admissible densities, defined as a subset of $H^{1}(\Omega)$,
given by
\begin{equation}
\mathcal{A}=\left\{ \rho\in H^{1}(\Omega):\delta_{\rho}\leq\rho\leq1\mbox{ a.e.}\right\} \label{eq:A}
\end{equation}
The proof of existence of minimizers for (\ref{eq:reg-compliance})
can be found in \citep{Talischi-JMPS} (see also \citep{Bendsoe-book}
for a weaker result) and essentially follows from compactness of the
minimizing sequences of (\ref{eq:reg-compliance}) in $L^{p}(\Omega)$,
$1\leq p<\infty$. We note that the norm of the density gradient also
appears in phase field formulations of topology optimization (see,
for example, \citep{Bourdin:2003p1799,Burger:2006p1864,Takezawa:2010p2730})
as an interfacial energy term and is accompanied by a double-well
potential penalizing intermediate densities. Taken together with appropriately
chosen coefficients, the two terms serve as approximation to the perimeter
of the design.

Under an additional assumption of $\partial\rho/\partial\mathbf{n}=0$
on $\partial\Omega$ and $\rho\in H^{2}(\Omega)$, the Tikhonov regularization
term can be written as $\frac{1}{2}\left\langle \rho,-\beta\Delta\rho\right\rangle $.
Similarly, the more general regularization term $\frac{1}{2}\left\langle \nabla\rho,\kappa\nabla\rho\right\rangle $
in which $\kappa(\mathbf{x})$ is a bounded and positive-definite
matrix prescribing varying regularity of $\rho$ in $\Omega$ can
be written as $\frac{1}{2}\left\langle \rho,-\nabla\cdot\left(\kappa\nabla\rho\right)\right\rangle $.
For brevity and emphasizing the quadratic form of this type of regularization,
in the next two sections, we write the regularizer generically as
\begin{equation}
\frac{1}{2}\left\langle \rho,\mathcal{R}\rho\right\rangle \label{eq:generic-reg}
\end{equation}
where $\mathcal{R}$ is a linear, self-adjoint and positive semi-definite
operator on $\mathcal{A}$, though the additional assumption on densities
are in fact not required. 

Finally, we recall that the gradient of compliance (with respect to
variations of density in the $L^{2}$-metric) is given by \citep{Bendsoe-book}
\begin{equation}
J'(\rho)=-E(\rho)+\lambda\label{eq:J_prime}
\end{equation}
where $E(\rho)=p\rho^{p-1}\mathbf{C}_{0}\boldsymbol{\epsilon}(\mathbf{u}_{\rho}):\boldsymbol{\boldsymbol{\epsilon}}(\mathbf{u}_{\rho})$
is a strain energy density field. Note that $E(\rho)$ is non-negative
for any admissible density and this is related to the monotonicity
of the self-adjoint compliance problem: given densities $\rho_{1}$
and $\rho_{2}$ such that $\rho_{1}\leq\rho_{2}$ a.e., one can show
$\ell(\mathbf{u}_{\rho_{1}})\geq\ell(\mathbf{u}_{\rho_{2}})$. This
property is the main reason why we restrict our attention in this
paper to compliance minimization (though in section 7, we will provide
an example of compliant mechanism design which is not self-adjoint).
Observe that $\hat{\rho}$ is a stationary point of $J$ if

\begin{equation}
\begin{cases}
E(\hat{\rho})(\mathbf{x)}<\lambda, & \mbox{if }\hat{\rho}(\mathbf{x)}=\delta_{\rho}\\
E(\hat{\rho})(\mathbf{x)}=\lambda, & \mbox{if }\delta_{\rho}<\hat{\rho}(\mathbf{x)}<1\\
E(\hat{\rho})(\mathbf{x)}>\lambda, & \mbox{if }\hat{\rho}(\mathbf{x)}=1
\end{cases}\label{eq:cond-opt}
\end{equation}
Thus, in regions where $E(\hat{\rho})$ exceeds the penalty parameter
$\lambda$ (regions that experience ``large'' deformation), density
is at its maximum. Similarly, below this cutoff value the density
is equal to the lower bound $\delta_{\rho}$. Everywhere else, i.e.,
in the regions of intermediate density, the strain energy density
is equal to the penalty parameter $\lambda$.

\begin{figure}
\centering{}\includegraphics[scale=0.62]{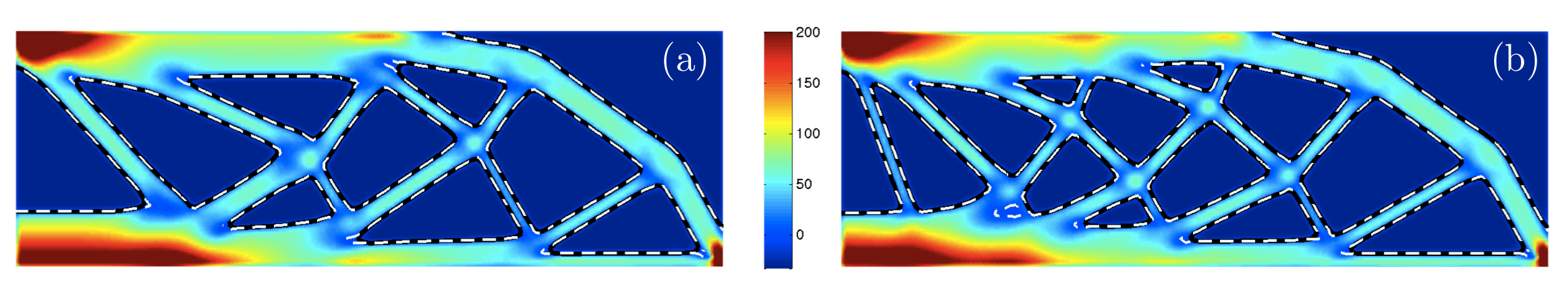}\caption{Plot of $E(\rho)-\lambda$ for two solutions to the MBB beam problem
with $\beta=0.01$ (a) corresponds to solution shown in Figure \ref{fig:Mbb-01}(b)
and (b) corresponds to solution shown in Figure \ref{fig:Mbb-01}(c).
The black line is the contour line for $\rho=1/2$ and the dashed
white line is the contour line where $E(\rho)=\lambda$. Note that
only half the design domain is shown and the range of the colorbar
is limited to $\left[-\lambda,6\lambda\right]$ for better visualization.\label{fig:StrainEnergy}}
\end{figure}

Figure \ref{fig:StrainEnergy} shows the distribution of $E(\rho)-\lambda$
for solutions to (\ref{eq:reg-compliance}) obtained using the proposed
algorithm (cf. section 7 and Figures \ref{fig:Mbb-01}(b) and (c)).
Superimposed are the contour lines associated with $\rho=1/2$ (plotted
in black) representing the boundary of the optimal shape and $E(\rho)=\lambda$
(plotted in dashed white). The fact that these lines are nearly coincident
shows that the solutions to the regularized problem, at least for
sufficiently small regularization parameter $\beta$, are close to
ideal in the sense that they nearly satisfy the stationarity condition
for the structural objective $J$.

\section{General Splitting Algorithm}

In this section, we discuss a generalization of the forward-backward
splitting algorithm that was explored in \citep{Talischi-JMPS} for
solving the regularized compliance minimization problem. The key idea
behind this and other similar decomposition methods \citep{Cohen:1978p3656,Chen:1997p3679,Patriksson:1998p3659}
is the separate treatment of constituent terms of the cost function. 

A general algorithm for finding a minimizer of $\tilde{J}(\rho)$
consists of subproblems of the form:
\begin{equation}
\rho_{n+1}=\underset{\rho\in\mathcal{A}_{n}}{\mbox{argmin}}\ J(\rho_{n})+\left\langle \rho-\rho_{n},J'(\rho_{n})\right\rangle +\frac{1}{2\tau_{n}}\left\Vert \rho-\rho_{n}\right\Vert _{\mathcal{H}_{n}}^{2}+\frac{1}{2}\left\langle \rho,\mathcal{R}\rho\right\rangle \label{eq:gen-subproblem}
\end{equation}
where $\mathcal{H}_{n}$ is a bounded and positive-definite linear
operator. Compared to (\ref{eq:reg-compliance}), we can see that
while the regularization term has remained intact, $J$ is replaced
by a local quadratic model around $\rho_{n}$ in which $\mathcal{H}_{n}$
may be viewed as an approximation to the Hessian of $J$ evaluated
at $\rho_{n}$. Note that constant terms such as $J(\rho_{n})$ and
$\left\langle \rho_{n},J'(\rho_{n})\right\rangle $ do not affect
the optimization but are provided to emphasize the expansion of $J$.
Moreover, $\tau_{n}>0$ is a \emph{step size} parameter that determines
the curvature of this approximation. For sufficiently small $\tau_{n}$
(large curvature), the approximation is conservative in that it majorizes
(lies above) $J$, which is crucial in guaranteeing decent in each
iteration and overall convergence of the algorithm (see section 6).

We have included another limiting measure in (\ref{eq:gen-subproblem}),
a minor departure from the above-mentioned references, by replacing
the constraint set $\mathcal{A}$ by a subset $\mathcal{A}_{n}$ in
order to limit the point-wise change in the density to a specified
\emph{move limit} $m_{n}$. More specifically, we have defined 
\begin{equation}
\mathcal{A}_{n}=\left\{ \rho\in\mathcal{A}:\left|\rho-\rho_{n}\right|\leq m_{n}\mbox{ a.e.}\right\} =\left\{ \rho\in H^{1}(\Omega):\rho_{n}^{\mathsf{L}}\leq\rho\leq\rho_{n}^{\mathsf{U}}\mbox{ a.e.}\right\} \label{eq:A_n}
\end{equation}
where in the latter expression
\begin{equation}
\rho_{n}^{\mathsf{L}}=\delta_{\rho}\wedge\left(\rho_{n}-m_{n}\right),\quad\rho_{n}^{\mathsf{U}}=1\vee\left(\rho_{n}+m_{n}\right)\label{eq:rho_n_L_U}
\end{equation}
The presence of move limits (akin to a trust region strategy) is common
in topology optimization literature as a means to stabilize the topology
optimization algorithm, especially in the early iterations to prevent
members from forming too prematurely. As we will show with an example,
this is only important when a smaller regularization parameter is
used and the final topology is complex. Near the optimal solution,
the move limit strategy is typically inoperative. Of course, by setting
$m_{n}\equiv1$, we can get $\mathcal{A}=\mathcal{A}_{n}$ and recover
the usual form of (\ref{eq:gen-subproblem}).

Ignoring the constant terms and with simple rearrangement, we can
show that (\ref{eq:gen-subproblem}) is equivalent to
\begin{equation}
\rho_{n+1}=\underset{\rho\in\mathcal{A}_{n}}{\mbox{argmin}}\ \left\Vert \rho-\rho_{n+1}^{*}\right\Vert _{\left(\mathcal{H}_{n}+\tau_{n}\mathcal{R}\right)}^{2}\label{eq:sub-alt-form1}
\end{equation}
where the \emph{interim density} $\rho_{n+1}^{*}$ is given by
\begin{equation}
\rho_{n+1}^{*}=\left(\mathcal{H}_{n}+\tau_{n}\mathcal{R}\right)^{-1}\left[\mathcal{H}_{n}\rho_{n}-\tau_{n}J'(\rho_{n})\right]\label{eq:interim-1}
\end{equation}
Alternatively, the interim density can be written as a Newton-type
update where the gradient of $\tilde{J}$ is scaled by the inverse
of its approximate Hessian, namely 
\begin{equation}
\rho_{n+1}^{*}=\rho_{n}-\tau_{n}\left(\mathcal{H}_{n}+\tau_{n}\mathcal{R}\right)^{-1}\left[J'(\rho_{n})+\mathcal{R}\rho_{n}\right]\label{eq:interim-2}
\end{equation}
Returning to (\ref{eq:sub-alt-form1}), we can see that next density
$\rho_{n+1}$ is defined as the \emph{projection} of the interim density,
with respect to the norm defined by $\mathcal{H}_{n}+\tau_{n}\mathcal{R}$,
onto the constraint space $\mathcal{A}_{n}$. From the assumptions
on properties of $\mathcal{H}_{n}$ and the Tikhonov regularization
operator $\mathcal{R}$ and the fact that $\mathcal{A}_{n}$ is a
closed convex subset of $H^{1}(\Omega)$, it follows that the projection
is well-defined and there is a unique update $\rho_{n+1}$. 

By setting $\mathcal{R}=-\beta\Delta$, which corresponds to the regularization
term of $\eqref{eq:reg-compliance}$ and choosing $\mathcal{H}_{n}$
to be the identity map $\mathcal{I}$, we recover the forward-backward
algorithm investigated in \citep{Talischi-JMPS}. In this case, the
interim update satisfies the Helmholtz equation
\begin{equation}
\left(\mathcal{I}-\tau_{n}\beta\Delta\right)\rho_{n+1}^{*}=\rho_{n}-\tau_{n}J'(\rho_{n})\label{eq:Helmh-filter}
\end{equation}
with homogenous Neumann boundary conditions. Note that the right hand
side is the usual gradient descent step (with step size $\tau_{n}$)
associated with $J$ (the forward step) and the interim density is
obtained from application of the inverse of the Helmholtz operator
(the backward step), which can be viewed as the filtering of right-hand-side
with the Gaussian Green's function of the Helmholtz equation%
\footnote{The designations ``forward'' and ``backward'' step come from the
fact that (\ref{eq:Helmh-filter}) can be written as $\rho_{n+1}^{*}=\left(\mathcal{I}+\tau_{n}\mathcal{R}\right)^{-1}\left(\mathcal{I}-\tau_{n}J'\right)\rho_{n}$.
Similarly, (\ref{eq:interim-1}) has equivalent expression $\rho_{n+1}^{*}=\left(\mathcal{I}+\tau_{n}\mathcal{H}_{n}^{-1}\mathcal{R}\right)^{-1}\left(\mathcal{I}-\tau_{n}\mathcal{H}_{n}^{-1}J'\right)\rho_{n}$.%
}. As mentioned in the introduction, this appearance of filtering is
fundamentally different from density and sensitivity filtering methods.
Moreover, the projection operation in this case is with respect to
a scaled Sobolev metric, namely
\begin{equation}
\rho_{n+1}=\underset{\rho\in\mathcal{A}_{n}}{\mbox{argmin}}\ \left\Vert \rho-\rho_{n+1}^{*}\right\Vert ^{2}+\beta\tau_{n}\left|\rho-\rho_{n+1}^{*}\right|_{1}^{2}\label{eq:scaled-sob-proj}
\end{equation}
which numerically requires the solution to a box-constrained convex
quadratic program. In \citep{Talischi-JMPS}, we also explored an
``inconsistent'' variation of this algorithm where we neglected
the second term in (\ref{eq:scaled-sob-proj}) and essentially used
the $L^{2}$-metric for the projection step. Due to the particular
geometry of the box constraints in $\mathcal{A}_{n}$, the $L^{2}$-projection
has the explicit solution given by 
\begin{equation}
\rho_{n+1}=\left(\rho_{n+1}^{*}\wedge\rho_{n}^{\mathsf{L}}\right)\vee\rho_{n}^{\mathsf{U}}\label{eq:L2-proj}
\end{equation}
The appeal of this min/max type operation is that it is trivial from
the computational point of view. Moreover, it coincides with the last
step in the OC update scheme \citep{Bendsoe-book}. However, this
is an inconsistent step for Tikhonov regularized problem since $\rho_{n+1}$
need not lie in $H^{1}(\Omega)$. In fact, strictly speaking, (\ref{eq:L2-proj})
is valid only if $\mathcal{A}_{n}$ is enlarged from functions in
$H^{1}(\Omega)$ to all functions in $L^{2}(\Omega)$ bounded below
by $\rho_{n}^{\mathsf{L}}$ and above by $\rho_{n}^{\mathsf{U}}$.
In spite of this inconsistency, the algorithm composed of (\ref{eq:Helmh-filter})
and (\ref{eq:L2-proj}) was convergent and numerically shown to produce
noteworthy solutions with minimal intermediate densities. This merits
a separate investigation since as suggested in \citep{Talischi-JMPS},
this algorithm may in fact solve a smoothed version of the perimeter
constraint problem where the regularization term is the total variation
of the density field. We will return to the use of $L^{2}$-projection
later in section 6 but \emph{this time in a consistent manner with
the aid of the two-metric projection approach of \citep{Bert-ProjNewton,Gafni-TMP}}.

\section{Optimality Criteria and Sensitivity Filtering}

In structural optimization, the optimality criteria (OC) method is
preferred to the gradient descent algorithm since it typically enjoys
faster convergence (see \citep{Arora:1980p3660} on the relationship
between the two methods). Our interest here in the OC method is that
the density and sensitivity filtering methods are typically implemented
in the OC framework. Moreover, as we shall see, this examination will
lead to the choice of $\mathcal{H}_{n}$ in the algorithm (\ref{eq:gen-subproblem}). 

The interim density in the OC method for the compliance minimization
problem (in the absence of regularization) is obtained from the fixed
point iteration
\begin{equation}
\rho_{n+1}^{*}=\rho_{n}\left[\frac{E(\rho_{n})}{\lambda}\right]^{1/2}\equiv\rho_{n}\left[e_{\lambda}(\rho_{n})\right]^{1/2}\label{eq:OC-update}
\end{equation}
Note that the strain energy density $E(\rho_{n})$ and subsequently
its normalization $e_{\lambda}(\rho_{n})$ are non-negative for any
admissible density $\rho_{n}$ and therefore $\rho_{n+1}^{*}$ is
well-defined. Recalling the necessary condition of optimality for
an optimal density $\hat{\rho}$ stated in (\ref{eq:cond-opt}), it
is evident that such $\hat{\rho}$ is a fixed point of the OC iteration.
Intuitively, the current density $\rho_{n}$ is increased (decreased)
in regions where $E(\rho_{n})$ is greater (less) than the penalty
parameter $\lambda$ by a factor of $\left[e_{\lambda}(\rho_{n})\right]^{1/2}$.
The next density $\rho_{n+1}$ in the OC is given by (\ref{eq:L2-proj}). 

It is more useful here to adopt an alternative view of the OC scheme,
namely that the OC update can be seen as the solution to an approximate
subproblem where compliance is replaced by a Taylor expansion in the
intermediate field $\rho^{-1}$ \citep{Groenwold:2008p2685}. The
intuition behind such expansion is that locally compliance is inversely
proportional to density. In particular, $\rho_{n+1}^{*}$ can be shown
to be the stationary point of the ``reciprocal approximation'' around
$\rho_{n}$ defined by

\begin{equation}
J_{\mathsf{rec}}(\rho;\rho_{n})\equiv\ell(\mathbf{u}_{\rho_{n}})+\left\langle \frac{\rho_{n}}{\rho}\left(\rho-\rho_{n}\right),-E(\rho_{n})\right\rangle +\lambda\int_{\Omega}\rho\mathrm{d}\mathbf{x}\label{eq:J_rec}
\end{equation}
Note that the expansion in the inverse of density is carried out only
for the compliance term, and the volume term, which is already linear,
is not altered. The expression for $J_{\mathsf{rec}}(\rho;\rho_{n})$
can be alternatively written as
\begin{equation}
J_{\mathsf{rec}}(\rho;\rho_{n})=J(\rho_{n})+\left\langle \rho-\rho_{n},J'(\rho_{n})\right\rangle +\frac{1}{2}\left\langle \rho-\rho_{n},\frac{2E(\rho_{n})}{\rho}\left(\rho-\rho_{n}\right)\right\rangle \label{eq:J_rec_alt}
\end{equation}
which highlights the fact that the (nonlinear) curvature term in (\ref{eq:J_rec_alt})
makes it a more accurate approximation of compliance compared to the
linear expansion. With regard to the OC update, one can show that
the interim update satisfies $J'_{\mathsf{rec}}(\rho_{n+1}^{*};\rho_{n})=0$,
and its $L^{2}$-projection is indeed the minimizer of $J_{\mathsf{rec}}(\rho;\rho_{n})$
over $\mathcal{A}_{n}$ (again enlarged to $L^{2}$).

\begin{figure}
\centering{}\includegraphics[scale=0.6]{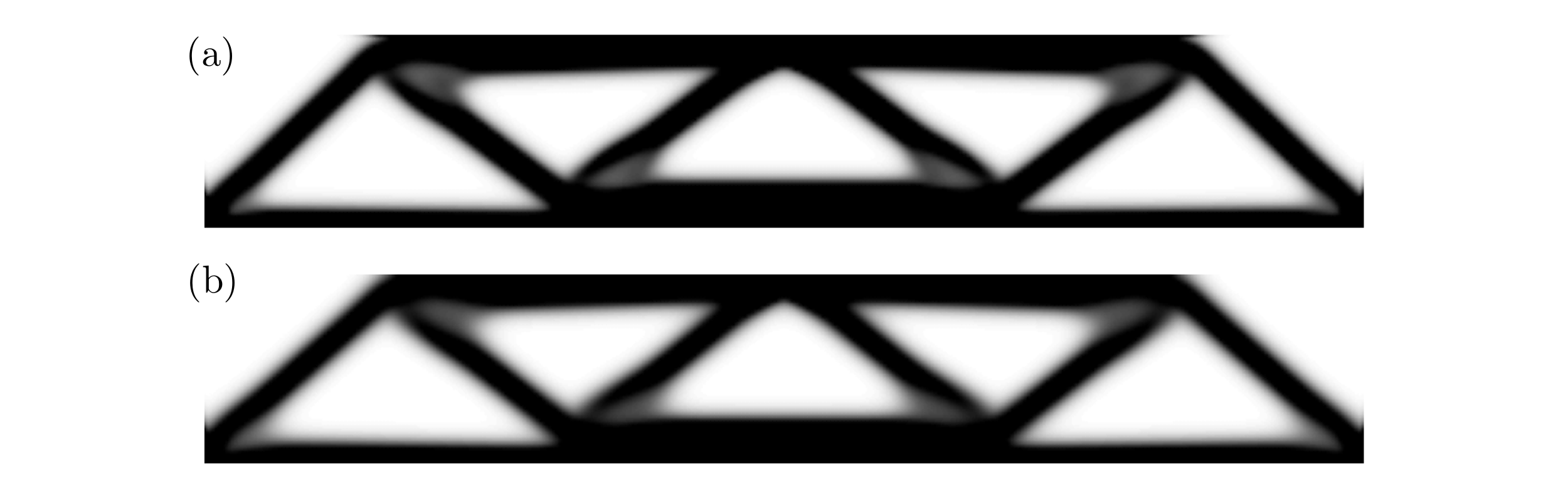}\caption{(a) The solution to the MBB beam problem (see section 6) using the
sensitivity filtering method (consisting of (\ref{eq:interim-sens-filt})
and (\ref{eq:L2-proj})) (b) The solution using the update steps (\ref{eq:FBS-OC-subs})
and (\ref{eq:L2-proj}). In both cases, $\mathcal{F}$ is taken to
be the ``Helmholtz'' filter and the move limit was set to $m_{n}=0.25$
\label{fig:ResultSenOCvsModFBS}}
\end{figure}

We now turn to the \emph{sensitivity filtering} method, which is described
with the OC algorithm. Let $\mathcal{F}$ denote a linear filtering
map, for example, the Helmholtz filter $\mathcal{F}=\left(\mathcal{I}-r^{2}\Delta\right)^{-1}$
discussed before or the convolution filter of radius radius $r$ \citep{Bourdin:2001p457,Borrvall:2001p422}
\begin{equation}
\mathcal{F}(\psi)(\mathbf{x})\equiv\int_{\Omega}F_{r}(\mathbf{x}-\mathbf{y})\psi(\mathbf{y})\mathrm{d}\mathbf{y}\label{eq:linear-filer}
\end{equation}
where the kernel is the linear hat function $F_{r}(\mathbf{x})=\max\left(1-\left|\mathbf{x}\right|/r,0\right)$.
The main idea in the sensitivity filtering method is that $e_{\lambda}(\rho_{n})$
is heuristically replaced by the following smoothed version%
\footnote{Notice that the filtering map is applied to the scaling of $e_{\lambda}(\rho_{n})$
by the density field itself, which is not easy to explain/justify.%
} 
\begin{equation}
\tilde{e}_{\lambda}(\rho_{n})\equiv\frac{1}{\rho_{n}}\mathcal{F}\left[\rho_{n}e_{\lambda}(\rho_{n})\right]\label{eq:E_tilde}
\end{equation}
before entering the OC update. The interim density update is thus
given by
\begin{equation}
\rho_{n+1}^{*}=\rho_{n}\left[\tilde{e}_{\lambda}(\rho_{n})\right]^{1/2}=\rho_{n}\left\{ \frac{\mathcal{F}\left[\rho_{n}e_{\lambda}(\rho_{n})\right]}{\rho_{n}}\right\} ^{1/2}=\rho_{n}^{1/2}\mathcal{F}\left[\rho_{n}e_{\lambda}(\rho_{n})\right]^{1/2}\label{eq:interim-sens-filt}
\end{equation}
A key observation in this work is that if \emph{we replace the gradient
decent step in forward-backward algorithm (cf. (\ref{eq:Helmh-filter}))
with the OC step, we obtain a similar update scheme to that of the
sensitivity filtering method}. More specifically, note that (\ref{eq:Helmh-filter})
can be written as $\rho_{n+1}^{*}=\mathcal{F}\left[\rho_{n}-\tau_{n}J'(\rho_{n})\right]$.
Substituting the term in the bracket with $\rho_{n}\left[e(\rho_{n})\right]^{1/2}$
gives
\begin{equation}
\rho_{n+1}^{*}=\mathcal{F}\left\{ \rho_{n}\left[e_{\lambda}(\rho_{n})\right]^{1/2}\right\} \label{eq:FBS-OC-subs}
\end{equation}
which resembles (\ref{eq:interim-sens-filt}). In fact, as illustrated
in Figure \ref{fig:ResultSenOCvsModFBS}, the two expressions produce
very similar final results (in particular, observe the similarity
between the patches of intermediate density in the corners that is
characteristic of the sensitivity filtering method). Of course, the
leap from the forward-backward algorithm to (\ref{eq:FBS-OC-subs}),
just like the sensitivity filtering method, lacks mathematical justification.
However, we will expand upon this observation and next derive the
algorithm similar to this empirical modification of the forward-backward
algorithm in a consistent manner. 

\subsection*{Embedding the Reciprocal Approximation}

Recalling the role of the reciprocal approximation of compliance in
the OC method, the key idea is to embed such an approximation in the
general subproblem of (\ref{eq:gen-subproblem}). We do so by choosing
$\mathcal{H}_{n}$ to be the Hessian of $J_{\mathsf{rec}}(\rho;\rho_{n})$
evaluated at $\rho_{n}$, namely%
\footnote{We note that the use of a quadratic approximations of the reciprocal
approximation has also been pursued in \citep{Groenwold:2009p2639,Groenwold:2010p2681}.%
}
\begin{equation}
\mathcal{H}_{n}=J_{\mathsf{rec}}''(\rho_{n};\rho_{n})=\frac{2E(\rho_{n})}{\rho_{n}}\mathcal{I}\label{eq:Hn_rec_app}
\end{equation}
As noted earlier, $E(\rho)$ is a non-negative function for any admissible
$\rho$ but may vanish in some subset of $\Omega$. This means that
$\mathcal{H}_{n}$ is only positive semi-definite and does not satisfy
the definiteness requirement for use in (\ref{eq:gen-subproblem}).
We can remedy this by replacing $E(\rho_{n})$ in (\ref{eq:Hn_rec_app})
with $E(\rho_{n})\wedge\delta_{E}$ where $0<\delta_{E}\ll\lambda$
is a prescribed constant. However, in most compliance problems (e.g.,
the benchmark problem considered later in section 7) the strain energy
field is strictly positive for all admissible densities. In fact,
the regions with zero strain energy density do not experience any
deformation and in light of the conditions of optimality (\ref{eq:cond-opt})
should be assigned the minimum density. Therefore, to simplify the
matters, \emph{we assume in the remainder of this section that the
loading and support conditions defined on $\Omega$ are such that
$E(\rho)\geq\delta_{E}$ almost everywhere} for all $\rho\in L^{\infty}(\Omega;\left[\delta_{\rho},1\right])$.

Comparing the quadratic approximation of $J$ with this choice of
$\mathcal{H}_{n}$ and the reciprocal approximation itself (cf. (\ref{eq:J_rec_alt})),
we see that the difference is in their curvature terms (the linear
terms of course match). The curvature of the quadratic model depends
on and can be controlled by $\tau_{n}$ while the nonlinear curvature
in $J_{\mathrm{rec}}$ is a function of $\rho$.

\begin{figure}
\centering{}\includegraphics[scale=0.5]{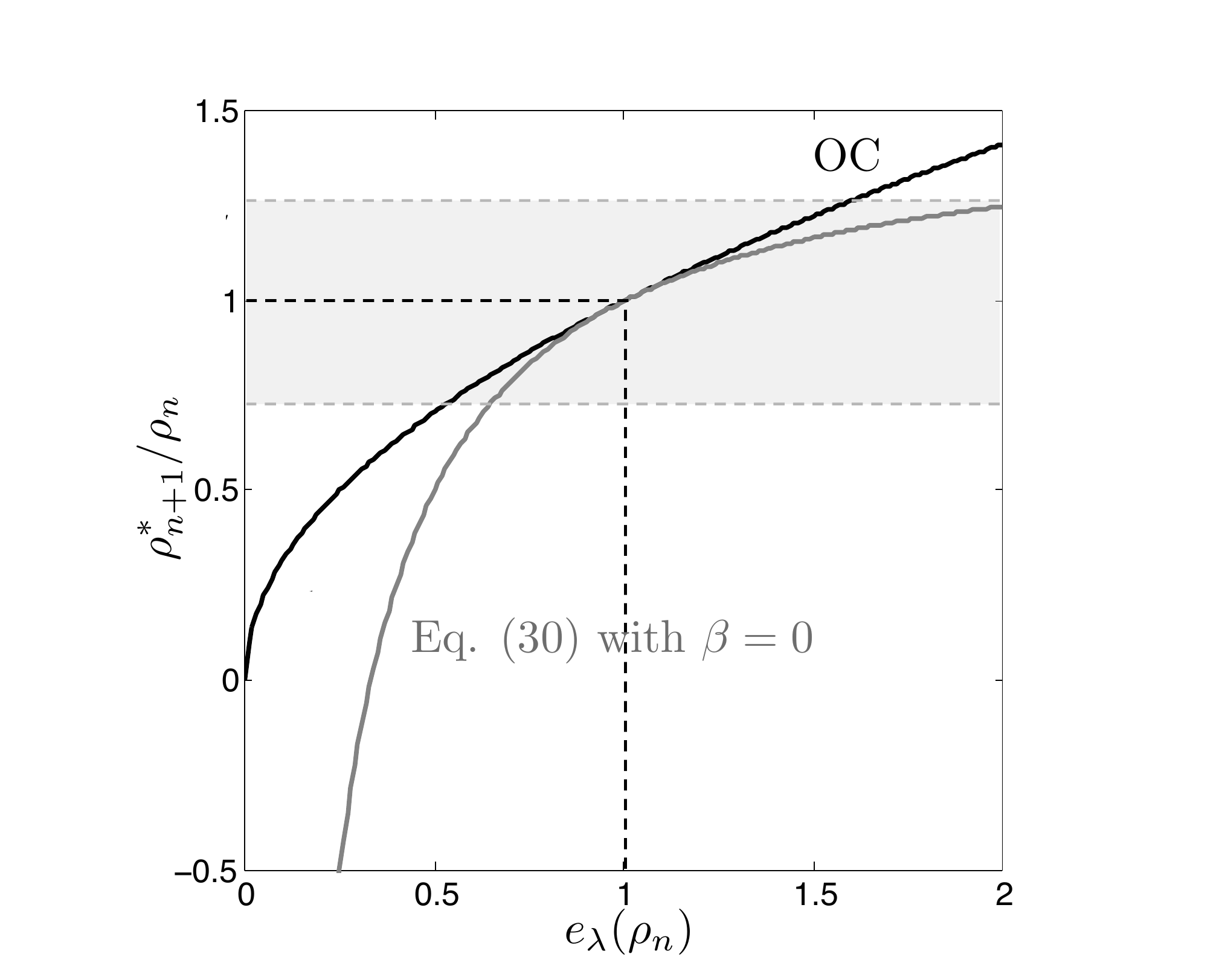}\caption{Comparison between scaling terms appearing in the OC update and right
hand side of (\ref{eq:subs-Hn-inter3}). The OC is more aggressive
in regions $e_{\lambda}(\rho_{n})>1$ and less aggressive when $e_{\lambda}(\rho_{n})<1$.
\label{fig:OC-update}}
\end{figure}

Substituting (\ref{eq:Hn_rec_app}) into (\ref{eq:interim-1}), the
expression for the interim density becomes
\begin{equation}
\left[\frac{2E(\rho_{n})}{\rho_{n}}\mathcal{I}+\tau_{n}\mathcal{R}\right]\rho_{n+1}^{*}=2E(\rho_{n})+\tau_{n}\left[E(\rho_{n})-\lambda\right]=\left(2+\tau_{n}\right)E(\rho_{n})-\tau_{n}\lambda\label{eq:subst-Hn-inter1}
\end{equation}
Multiplying by $\rho_{n}/\left[2E(\rho_{n})\right]$ and simplifying
yields 
\begin{equation}
\left[\mathcal{I}+\frac{\rho_{n}}{2E(\rho_{n})}\tau_{n}\mathcal{R}\right]\rho_{n+1}^{*}=\rho_{n}\left[\left(1+\frac{\tau_{n}}{2}\right)-\frac{\tau_{n}}{2e_{\lambda}(\rho_{n})}\right]\label{eq:subst-Hn-inter2}
\end{equation}
To better understand the characteristics of this update, let us specialize
to the case of Tikhonov regularization and set $\tau_{n}=1$ (so that
the quadratic model and the reciprocal approximation have the same
curvature at $\rho_{n}$). This gives 
\begin{equation}
\left[\mathcal{I}-\frac{\rho_{n}}{2E(\rho_{n})}\beta\Delta\right]\rho_{n+1}^{*}=\rho_{n}\left[\frac{3}{2}-\frac{1}{2e_{\lambda}(\rho_{n})}\right]\label{eq:subs-Hn-inter3}
\end{equation}
First note that in the absence of regularization (i.e., $\beta=0$),
the update relation has the same fixed-point iteration form as the
OC update with the ratio $e_{\lambda}(\rho_{n})$ determining the
scaling of $\rho_{n}$. The scaling field here is $3/2-1/\left[2e_{\lambda}(\rho_{n})\right]$
whereas in the OC method it is given by $\left[e_{\lambda}(\rho_{n})\right]^{1/2}$.
As shown in Figure \ref{fig:OC-update}, the scaling fields and their
derivatives coincide in the regions where $e_{\lambda}(\rho_{n})=1$,
which means that locally the two are similar. The reduction in density
is more aggressive with this scaling when $e_{\lambda}(\rho_{n})<1$
whereas the OC update leads to larger increase for $e_{\lambda}(\rho_{n})>1$.

As with the forward-backward algorithm (cf. equation (\ref{eq:Helmh-filter})),
the presence of regularization again leads to the appearance of Helmholtz
filtering (the inverse of left-hand-side operator) but with two notable
differences. First, the right-hand-side term now is an OC-like scaling
of density instead of the gradient descent step (the same is true
in (\ref{eq:subst-Hn-inter2}) for an arbitrary step size $\tau_{n}$).
Furthermore, the filtering is not uniform across the domain and its
degree of smoothening is scaled by $\rho_{n}/\left[2E(\rho_{n})\right]$.
The important result here is that, by embedding the reciprocal approximation
of compliance in our quadratic model, we are able to obtain a relation
for the $\rho_{n+1}^{*}$ that features an OC-like right-hand-side
and its filtering, very much similar in form to the (heuristically)
fabricated update scheme of (\ref{eq:FBS-OC-subs}) that was compared
to the sensitivity filtering. 

Another key difference between the forward-backward algorithm and
the OC-based filtering methods is that the projection of $\rho_{n+1}^{*}$
defining the next iterate $\rho_{n+1}$ in the forward-backward scheme
is with respect to the metric induced by $\mathcal{H}_{n}+\tau_{n}\mathcal{R}$
in contrast to the $L^{2}$-projection given by (\ref{eq:L2-proj}).
As discussed before, the $L^{2}$-projection is well-suited for the
geometry of the constraint set $\mathcal{A}_{n}$ due to decomposition
of box constraints. It may be tempting to inconsistently use the interim
density (\ref{eq:subst-Hn-inter2}) with the $L^{2}$-projection but
this is not necessarily guaranteed to decrease the cost function%
\footnote{Numerically one would observe that such an inconsistent algorithm
excessively removes material and leads to final solutions with low
volume fraction%
}. Arbitrary projections of unconstrained Newton steps is not mathematically
warranted.

In section 6, we explore a variant of the splitting algorithm that
is related to the two-metric projection method of \citep{Bert-ProjNewton,Gafni-TMP},
and allows for the use of a more convenient metric for the projection
step. This can be done provided that the operator whose associated
norm defines the gradient%
\footnote{Recall that $\mathcal{B}^{-1}f'(\rho)$ is the gradient of functional
$f$ with respect to the metric induced by $\mathcal{B}$. As such,
Newton's method and its variations (such as the present framework)
can be thought of as gradient descent algorithms with respect to a
variable metric defined by the (approximate) Hessian.%
} is modified appropriately in the regions where the constraints are
active. More specifically, in the interim update step (cf. (\ref{eq:interim-2})),
$\mathcal{H}_{n}+\tau_{n}\mathcal{R}$ is modified to produce a linear
operator $\mathcal{D}_{n}$ with a particular structure that eliminates
the coupling between regions of active and free constraints. The projection
of the interim density given by 
\begin{equation}
\rho_{n+1}^{*}=\rho_{n}-\tau_{n}\mathcal{D}_{n}^{-1}\left[J'(\rho_{n})+\mathcal{R}\rho_{n}\right]\label{eq:proj-Newt-cont}
\end{equation}
with respect to the $L^{2}$-norm is then guaranteed to decrease the
cost function%
\footnote{There is the technical issue that $L^{2}$-projection on a subset
of $H^{1}(\Omega)$ is not well-defined, which is why we defer the
exact outline of the two-metric projection method to the discrete
setting where this issue does not arise.%
}. Note that when there are no active constraints (e.g., in the beginning
of the algorithm the density field takes mostly intermediate values),
$\mathcal{D}_{n}=\mathcal{H}_{n}+\tau_{n}\mathcal{R}$ and (\ref{eq:subst-Hn-inter2})
holds for the interim update and its $L^{2}$-projection produces
the next iterate. In general, (\ref{eq:subst-Hn-inter2}) holds locally
for the regions where the box constraints are not active (i.e., regions
of intermediate density) and so the analogy to the sensitivity filtering
method holds in such regions. 

To avoid some technical nuisances (that the $L^{2}$-projection on
$H^{1}(\Omega)$ is not well-defined) and avoid the cumbersome notation
required to precisely define $\mathcal{D}_{n}$ in the continuum setting
(that may obscure the simple procedure for its construction), we defer
the details to section 6 where we describe the algorithm for the finite-dimensional
optimization problem obtained from the usual finite element approximation
procedure. The intuition developed in the preceding discussion carries
over to the discrete setting.

\section{Finite Element Approximation}

We begin with describing the approximate ``finite element'' optimization
problem, based on a typical choice of discretization spaces, and establish
the convergence of the corresponding optimal solutions to a solution
of the continuum problem (\ref{eq:reg-compliance}) in the limit of
mesh refinement. Our result proves strong convergence of a subsequence
of solutions, and therefore rules out the possibility of numerical
instabilities such as checkerboard patterns observed in density-based
methods. We remark that similar results are available for the density-based
restriction formulations (see for example \citep{Petersson:1998p440,Petersson:1999p452,Bourdin:2001p457})
and the proof is along the same lines. Such convergence results are
essential in justifying an overall optimization approach where one
first discretizes a well-posed continuum problem and then chooses
an algorithm to solve the resulting finite dimensional problem (this
is the procedure adopted in this work). Then, with the FE convergence
result in hand, the only remaining task is to analyze the convergence
of the proposed optimization algorithm, which is discussed in the
section 6.

\subsection{Convergence under mesh refinement}

Consider partitioning of $\Omega$ into pairwise disjoint finite elements
$\mathcal{T}_{h}=\left\{ \Omega_{e}\right\} _{e=1}^{l}$ with characteristic
mesh size $h$. Let $\mathcal{A}_{h}$ be the FE subspace of $\mathcal{A}$
based on this partition:
\begin{equation}
\mathcal{A}_{h}=\left\{ \rho\in C^{0}(\overline{\Omega}):\rho|_{\Omega_{e}}\in\mathcal{P}(\Omega_{e}),\forall e=1,\dots,l\right\} \cap\mathcal{A}\label{eq:A_h}
\end{equation}
where $\mathcal{P}(\Omega_{e})$ is a space of polynomial (rational
in the case of polygonal elements) functions defined on $\Omega_{e}$.
Similarly, we define:
\begin{equation}
\mathcal{V}_{h}=\left\{ \mathbf{u}\in C^{0}(\overline{\Omega};\mathbb{R}^{d}):\left[\mathbf{u}\right]_{i}|_{\Omega_{e}}\in\mathcal{P}(\Omega_{e}),\forall e=1,\dots,l,\forall i=1,\dots,d\right\} \cap\mathcal{V}\label{eq:V_h}
\end{equation}
We also assume that the mesh $\mathcal{T}_{h}$ is chosen in such
a way that the transition from $\Gamma_{D}$ to $\Gamma_{N}$ is properly
aligned with the mesh. In practice, both density and displacement
fields are discretized with linear elements (e.g., linear triangles,
bilinear quads or linearly-complete convex polygons in two spatial
dimensions). To avoid any ambiguity regarding the definition of the
FE partitions, we assume a regular refinement of the meshes such that
the resulting finite element spaces are ordered, e.g., $\mathcal{A}_{h}\supseteq\mathcal{A}_{h'}$
whenever $h\leq h'$. We consider the limit $h\rightarrow0$ to establish
convergence of solutions under mesh refinement. 

What is needed in the proof of convergence is the existence of an
interpolation operator $\mathcal{I}_{h}:\mathcal{V}\rightarrow\mathcal{V}_{h}$
such that for all $\mathbf{u}\in\mathcal{V}\cap H^{2}(\Omega;\mathbb{R}^{d})$
\begin{equation}
\left\Vert \mathcal{I}_{h}\mathbf{u}-\mathbf{u}\right\Vert _{1}\leq Ch\left|\mathbf{u}\right|_{2}\label{eq:interpolation-op}
\end{equation}
which in turn shows that $\mathcal{I}_{h}\mathbf{u}\rightarrow\mathbf{u}$
as $h\rightarrow0$. Similarly, we need the mapping $i_{h}:\mathcal{A}\rightarrow\mathcal{A}_{h}$
for the design space such that $i_{h}\rho\rightarrow\rho$ as $h\rightarrow0$.
The construction of such interpolants is standard in finite element
approximation theory, see for example \citep{Brenner-FE}.

The approximate finite element problem, specialized to Tikhonov regularization,
is defined by
\begin{equation}
\min_{\rho\in\mathcal{A}_{h}}\tilde{J}_{h}(\rho):=J_{h}(\rho)+\frac{\beta}{2}\left|\rho\right|_{1}^{2}\label{eq:min-Jt_Ah}
\end{equation}
where $J_{h}(\rho):=\ell(\mathbf{u}_{\rho,h})+\lambda\int_{\Omega}\rho\mathrm{d}\mathbf{x}$
and $\mathbf{u}_{\rho,h}$ is the solution to the Galerkin approximation
of (\ref{eq:weak-form}) given by
\begin{equation}
a(\mathbf{u}_{h},\mathbf{v}_{h};\rho)=\ell(\mathbf{v}_{h}),\qquad\forall\mathbf{v}_{h}\in\mathcal{V}_{h}\label{eq:FE-state-equation}
\end{equation}
By the principle of minimum potential, we can write 
\begin{equation}
\ell(\mathbf{u}_{\rho,h})=-2\min_{\mathbf{v}_{h}\in\mathcal{V}_{h}}\left[\frac{1}{2}a(\mathbf{v}_{h},\mathbf{v}_{h};\rho)-\ell(\mathbf{v}_{h})\right]=\max_{\mathbf{v}_{h}\in\mathcal{V}_{h}}\left[2\ell(\mathbf{v}_{h})-a(\mathbf{v}_{h},\mathbf{v}_{h};\rho)\right]\label{eq:compli-fe-energy-form}
\end{equation}
From the above relation, it is easy to see that $\mathcal{V}_{h}\subseteq\mathcal{V}$
implies $\ell(\mathbf{u}_{\rho,h})\leq\ell(\mathbf{u}_{\rho})$ for
any given $\rho$, and therefore
\begin{equation}
\tilde{J}_{h}(\rho)\leq\tilde{J}(\rho)\label{eq:Jh-compared-J}
\end{equation}
that is, the finite approximation of the state equation leads to a
smaller computed value of the cost function for any density field.

Consider a sequence of FE partitions $\mathcal{T}_{h}$ with $h\rightarrow0$
and let $\rho_{h}$ be the optimal solution to the associated finite
element approximation (\ref{eq:min-Jt_Ah}), i.e., minimizer of $\tilde{J}_{h}$
in $\mathcal{A}_{h}$. We first show the sequence $\rho_{h}$ is bounded
in $H^{1}(\Omega)$. To see this, fix $h_{0}$ in this sequence. If
$\hat{\rho}_{h}$ is the minimizer of $\tilde{J}$ in $\mathcal{A}_{h}$
(there is no approximation of the displacement field involved here),
then 
\begin{equation}
\tilde{J}(\hat{\rho}_{h})\leq\tilde{J}(\rho_{h_{0}})\label{eq:FE-conv-int-result1}
\end{equation}
since $\rho_{h_{0}}\in\mathcal{A}_{h_{0}}\subseteq\mathcal{A}_{h}$.
Now, from the definition of $\rho_{h}$ and (\ref{eq:Jh-compared-J}),
we have $\tilde{J}_{h}(\rho_{h})\leq\tilde{J}_{h}(\hat{\rho}_{h})\leq\tilde{J}(\hat{\rho}_{h})$
and so 
\begin{equation}
\tilde{J}_{h}(\rho_{h})\leq\tilde{J}(\rho_{h_{0}})=\tilde{J}_{h_{0}}(\rho_{h_{0}})+\left[J(\rho_{h_{0}})-J_{h_{0}}(\rho_{h_{0}})\right]:=\tilde{J}_{h_{0}}(\rho_{h_{0}})+\epsilon_{h_{0}}\label{eq:FE-conv-int-result2}
\end{equation}
where $\epsilon_{h_{0}}$ is the finite element error in computing
compliance of $\rho_{h_{0}}$ on mesh $\mathcal{T}_{h_{0}}$. Since
(\ref{eq:FE-conv-int-result2}) holds for all $h\leq h_{0}$, we conclude
that
\begin{equation}
\limsup_{h\rightarrow0}\tilde{J}_{h}(\rho_{h})\leq\tilde{J}_{h_{0}}(\rho_{h_{0}})+\epsilon_{h_{0}}\label{eq:FE-conv-int-result3}
\end{equation}
Both the compliance and volume terms in $\tilde{J}_{h}(\rho_{h})$
are uniformly bounded, and so (\ref{eq:FE-conv-int-result3}) shows
$\limsup_{h}\left|\rho_{h}\right|_{1}^{2}<\infty$ . Thus the sequence
$\rho_{h}$ is bounded in $H^{1}(\Omega)$. By Rellich's theorem \citet{Evans},
we have convergence of a subsequence, again denoted by $\left\{ \rho_{h}\right\} $,
strongly in $L^{2}(\Omega)$ and weakly in $H^{1}(\Omega)$ to some
$\rho^{*}\in\mathcal{A}$%
\footnote{To see that $\rho^{*}$ satisfies the bound constraints, we can consider
another subsequence for which the convergence is pointwise. %
}. We next show that $\rho^{*}$ is a solution to continuum problem,
thereby establishing the convergence of the FE approximate problem.
First note that by lower semi-continuity of the norm under weak convergence,
\begin{equation}
\left|\rho^{*}\right|_{1}^{2}\leq\liminf_{h}\left|\rho_{h}\right|_{1}^{2}\label{eq:lsc-semi-norm-weak}
\end{equation}
Furthermore, to show convergence of $\mathbf{u}_{\rho_{h},h}$ to
$\mathbf{u}_{\rho^{*}}$ in $H^{1}(\Omega;\mathbb{R}^{d})$, first
note that the convergence results stated in section 2 implies that
up to a subsequence $\mathbf{u}_{\rho_{h}}\rightarrow\mathbf{u}_{\rho^{*}}$
as $h\rightarrow0$. Moreover, 
\begin{eqnarray}
\left\Vert \mathbf{u}_{\rho^{*}}-\mathbf{u}_{\rho_{h},h}\right\Vert _{1} & \leq & \left\Vert \mathbf{u}_{\rho^{*}}-\mathbf{u}_{\rho_{h}}\right\Vert _{1}+\left\Vert \mathbf{u}_{\rho_{h}}-\mathbf{u}_{\rho_{h},h}\right\Vert _{1}\nonumber \\
 & \leq & \left\Vert \mathbf{u}_{\rho^{*}}-\mathbf{u}_{\rho_{h}}\right\Vert _{1}+\frac{M}{c}\left\Vert \mathbf{u}_{\rho_{h}}-\mathcal{I}_{h}(\mathbf{u}_{\rho_{h}})\right\Vert _{1}\label{eq:CeasLemma}\\
 & \leq & \left\Vert \mathbf{u}_{\rho^{*}}-\mathbf{u}_{\rho_{h}}\right\Vert _{1}+\hat{C}h\left|\mathbf{u}_{\rho_{h}}\right|_{2}\nonumber 
\end{eqnarray}
where the second inequality follows from Cea's lemma \citep{Brenner-FE}
and last inequality follows from estimate (\ref{eq:interpolation-op}).
Hence $\mathbf{u}_{\rho_{h},h}\rightarrow\mathbf{u}_{\rho^{*}}$ in
$H^{1}(\Omega;\mathbb{R}^{d})$ and so $J_{h}(\rho_{h})\rightarrow J(\rho^{*})$.
Together with the above inequality, we have
\begin{equation}
\tilde{J}(\rho^{*})\leq\liminf_{h}\tilde{J}_{h}(\rho_{h})\label{eq:fe-cont-opt-of-rhostar}
\end{equation}
To establish optimality of $\rho^{*}$, take any $\rho\in\mathcal{A}_{h}$.
The definition of $\rho_{h}$ as the optimal solution to (\ref{eq:min-Jt_Ah})
implies 
\begin{equation}
\tilde{J_{h}}(\rho_{h})\leq\tilde{J_{h}}\left[i_{h}(\rho)\right]\label{eq:fe-conv-int-result4}
\end{equation}
Using a similar argument as above, we can pass (\ref{eq:fe-conv-int-result4})
to the limit to show $\tilde{J}(\rho^{*})\leq\tilde{J}(\rho)$.

\subsection{The Discrete Problem}

We proceed to obtain explicit expressions for the discrete problem
(\ref{eq:min-Jt_Ah}) for a given finite element partition $\mathcal{T}_{h}$.
For each $\rho_{h}\in\mathcal{A}_{h}$, we have the expansion $\rho_{h}(\mathbf{x})=\sum_{k=1}^{m}\left[\mathbf{z}\right]_{k}\varphi_{k}(\mathbf{x})$
where $\mathbf{z}$ is the vector of nodal densities characterizing
$\rho_{h}$ and $\left\{ \varphi_{k}\right\} _{k=1}^{m}$ the set
of finite element basis functions for $\mathcal{A}_{h}$%
\footnote{Naturally we assume that the basis functions are such that for any
$\mathbf{z}\in\left[\delta_{\rho},1\right]^{m}$, the associated density
field lies in $\left[\delta_{\rho},1\right]$ everywhere. This is
satisfies, for example, if $0\leq\varphi_{k}\leq1$ for all $k$,
which is the case for linear convex $n$-gons \citep{Talischi:2010p3179}.%
}. The finite-dimensional space corresponding to $\mathcal{A}_{h}$
is simply the closed cube $\left[\delta_{\rho},1\right]^{m}$. Moreover,
the vector form for the Tikhonov regularization term is 
\begin{equation}
\frac{\beta}{2}\left|\rho_{h}\right|_{1}^{2}=\frac{1}{2}\mathbf{z}^{T}\mathbf{G}\mathbf{z}\label{eq:tikhonov-matrix-form}
\end{equation}
where $\mathbf{G}$ is the usual finite element matrix defined by
$\left[\mathbf{G}\right]_{k\ell}=\beta\int_{\Omega}\nabla\varphi_{k}\cdot\nabla\varphi_{\ell}\mathrm{d}\mathbf{x}$,
which is positive semi-definite. Similarly, the volume term $\int_{\Omega}\rho\mathrm{d}\mathbf{x}$
can be written as $\mathbf{z}^{T}\mathbf{v}$ where $\left[\mathbf{v}\right]_{k}=\int_{\Omega}\varphi_{k}\mathrm{d}\mathbf{x}$.

With regard to state equation (\ref{eq:FE-state-equation}), we make
one approximation in the energy bilinear form%
\footnote{This is a departure from the previous section but it can be accounted
for in the convergence analysis.%
} by assuming that the density field has a constant value over each
element, equal to the centroidal value, in the bilinear form. If $\mathbf{x}_{e}$
denotes the location of the centroid of element $\Omega_{e}$, we
replace each $\rho_{h}(\mathbf{x})$ by%
\footnote{Here $\chi_{A}$ is the characteristic function associated with set
$A$, i.e., a function that takes value of 1 for $\mathbf{x}\in A$
and zero otherwise.%
} 
\begin{equation}
\sum_{e=1}^{l}\chi_{\Omega_{e}}(\mathbf{x})\rho_{h}(\mathbf{x}_{e})\label{eq:element-wise-approx-dens}
\end{equation}
in the state equation. The use of piecewise element density is common
practice in topology optimization (cf. \citep{PolyTop}) and makes
the calculations and notation simpler. If $\left\{ \mathbf{N}_{i}\right\} _{i=1}^{q}$
denotes the basis functions for the displacement field such that $\mathbf{u}_{h}(\mathbf{x})=\sum_{i=1}^{q}\mathbf{\left[U\right]}_{i}\mathbf{N}_{i}(\mathbf{x})$,
the vector form of (\ref{eq:FE-state-equation}) is given by

\begin{equation}
\mathbf{K}\mathbf{U}=\mathbf{F}\label{eq:KU=00003DF}
\end{equation}
where the load vector $\left[\mathbf{F}\right]_{i}=\int_{\Gamma_{N}}\mathbf{t}\cdot\mathbf{N}_{i}\mathrm{d}s$
and the stiffness matrix, with the above approximation of density,
is
\begin{equation}
\left[\mathbf{K}\right]_{ij}=\int_{\Omega}\rho_{h}^{p}\mathbf{C}_{0}\nabla\mathbf{N}_{i}:\nabla\mathbf{N}_{j}\mathrm{d}\mathbf{x}=\sum_{e=1}^{l}\left[\rho_{h}(\mathbf{x}_{e})\right]^{p}\int_{\Omega_{e}}\mathbf{C}_{0}\nabla\mathbf{N}_{i}:\nabla\mathbf{N}_{j}\mathrm{d}\mathbf{x}\label{eq:Stiffness-matrix1}
\end{equation}
Let us define the matrix $\mathbf{P}$ whose $\left(e,k\right)$-entry
is given by $\mathbf{\left[P\right]}_{ek}=\varphi_{k}(\mathbf{x}_{e})$.
Then 
\begin{equation}
\rho_{h}(\mathbf{x}_{e})=\sum_{k=1}^{m}\varphi_{k}(\mathbf{x}_{e})\left[\mathbf{z}\right]_{k}=\sum_{k=1}^{m}\mathbf{\left[P\right]}_{ek}\left[\mathbf{z}\right]_{k}=\left[\mathbf{P}\mathbf{z}\right]_{e}\label{eq:Def-P}
\end{equation}
The vector $\mathbf{P}\mathbf{z}$ thus gives the vector of elemental
density values. Returning to (\ref{eq:Stiffness-matrix1}) and denoting
the \emph{element stiffness} matrix by $\mathbf{k}_{e}=\int_{\Omega_{e}}\mathbf{C}_{0}\nabla\mathbf{N}_{i}:\nabla\mathbf{N}_{j}\mathrm{d}\mathbf{x}$,
we have the simplified expression for the global stiffness matrix
\begin{equation}
\mathbf{K}(\mathbf{z})=\sum_{e=1}^{l}\left(\left[\mathbf{P}\mathbf{z}\right]_{e}\right)^{p}\mathbf{k}_{e}\label{eq:K-simple}
\end{equation}
The summation effectively represents the assembly routine in practice.
We note the continuity and ellipticity of the bilinear form (cf. (\ref{eq:coer-cont-of-a}))
and non-degeneracy of the finite element partition imply that the
eigenvalues of $\mathbf{K}(\mathbf{z})$ are bounded below by $c_{h}$
and above by $M_{h}$ (which depend on the mesh size -- see chapter
9 of \citep{Brenner-FE}) for all admissible density vectors $\mathbf{z}\in\left[\delta_{\rho},1\right]^{m}$. 

The discrete optimization problem (\ref{eq:min-Jt_Ah}) can now be
equivalently written as (with a slight abuse of notation for $J$
and $\tilde{J}$ ) 
\begin{equation}
\min_{\mathbf{z}\in\left[\delta_{\rho},1\right]^{m}}\tilde{J}(\mathbf{z}):=J(\mathbf{z})+\frac{1}{2}\mathbf{z}^{T}\mathbf{G}\mathbf{z}\label{eq:Discrete-opt}
\end{equation}
where 
\begin{equation}
J(\mathbf{z})=\mathbf{F}^{T}\mathbf{U}(\mathbf{z})+\lambda\mathbf{z}^{T}\mathbf{v}\label{eq:J-z}
\end{equation}
and $\mathbf{U(\mathbf{z})}$ is the solution to $\mathbf{K}(\mathbf{z})\mathbf{U}=\mathbf{F}$.
Observe that matrices $\mathbf{P}$ and $\mathbf{G}$, the vector
$\mathbf{v}$, as well as the element stiffness matrices $\mathbf{k}_{e}$
and load vector $\mathbf{F}$ are all fixed and do not change in the
course of optimization. Thus they can be computed once in the beginning
and stored.

The gradient of $J$ with respect to the nodal densities $\mathbf{z}$
can readily computed as
\begin{equation}
\partial_{k}J(\mathbf{z})=-\mathbf{U}(\mathbf{z})^{T}\left(\partial_{k}\mathbf{K}\right)\mathbf{U(\mathbf{z})}+\lambda\mathbf{\left[v\right]}_{k}\label{eq:dJ_k}
\end{equation}
The expression for $\partial_{k}\mathbf{K}$ can be obtained from
(\ref{eq:K-simple}). Defining the vector of strain energy densities
$\left[\mathbf{E}(\mathbf{z})\right]_{e}=p\left[\mathbf{P}\mathbf{z}\right]_{e}^{p-1}\mathbf{U}(\mathbf{z})^{T}\mathbf{k}_{e}\mathbf{U}(\mathbf{z})$,
we have
\begin{equation}
\nabla J(\mathbf{z})=-\mathbf{P}^{T}\mathbf{E}(\mathbf{z})+\lambda\mathbf{v}\label{eq:Grad_J}
\end{equation}
With the first order gradient information in hand, we can find the
reciprocal approximation%
\footnote{The reciprocal approximation to $f(\mathbf{x})$ at point $\mathbf{y}$
is given by $f(\mathbf{y})+\sum_{k=1}^{m}\left[x_{k}^{-1}y_{k}\left(x_{k}-y_{k}\right)\partial_{k}f(\mathbf{y})\right]$%
} of compliance about point $\mathbf{y}$ as
\begin{equation}
J_{\mathsf{rec}}(\mathbf{z};\mathbf{y})\equiv J(\mathbf{y})+\lambda\left(\mathbf{z}-\mathbf{y}\right)^{T}\mathbf{v}+\sum_{k=1}^{m}\left(\frac{\left[\mathbf{y}\right]_{k}}{\left[\mathbf{z}\right]_{k}}\right)\left(\left[\mathbf{z}\right]_{k}-\left[\mathbf{y}\right]_{k}\right)\left[-\mathbf{P}^{T}\mathbf{E}(\mathbf{y})\right]_{k}\label{eq:J_rec_z_y}
\end{equation}
The Hessian of $J_{\mathsf{rec}}(\mathbf{z};\mathbf{y})$, evaluated
at $\mathbf{z}=\mathbf{y}$, is a diagonal matrix with entries
\begin{equation}
h_{k}(\mathbf{y})=\partial_{kk}J_{\mathsf{rec}}(\mathbf{y};\mathbf{y})=\frac{2}{\left[\mathbf{y}\right]_{k}}\left[\mathbf{P}^{T}\mathbf{E}(\mathbf{y})\right]_{k},\quad k=1,\dots,m\label{eq:diag-rec-Hess}
\end{equation}
The entries of the vector $\mathbf{E}(\mathbf{y})$ are non-negative
for all admissible nodal densities but can be zero and therefore Hessian
of $J_{\mathsf{rec}}(\mathbf{z};\mathbf{y})$ is only positive semi-definite.

\section{Algorithms for the Discrete Problem}

We begin with the generalization of the forward-backward algorithm
for solving the discrete problem (\ref{eq:Discrete-opt}) before discussing
the two-metric projection variation. As in section 3, we consider
a splitting algorithm with iterations of the form
\begin{equation}
\mathbf{z}_{n+1}=\underset{\mathbf{z}_{n}^{\mathsf{L}}\leq\mathbf{z}\leq\mathbf{z}_{n}^{\mathsf{U}}}{\mbox{argmin}}\ Q_{J}(\mathbf{z};\mathbf{z}_{n},\tau_{n})+\frac{1}{2}\mathbf{z}^{T}\mathbf{G}\mathbf{z}\label{eq:z-Quad-Subproblem}
\end{equation}
where, compared to (\ref{eq:Discrete-opt}), the regularization term
is unchanged while $J$ is replaced by the following local quadratic
model around current iterate $\mathbf{z}_{n}$ 
\begin{equation}
Q_{J}(\mathbf{z};\mathbf{z}_{n},\tau_{n})=J(\mathbf{z}_{n})+\left(\mathbf{z}-\mathbf{z}_{n}\right)^{T}\nabla J(\mathbf{z}_{n})+\frac{1}{2\tau_{n}}\left\Vert \mathbf{z}-\mathbf{z}_{n}\right\Vert _{\mathbf{H}_{n}}^{2}\label{eq:QJ}
\end{equation}
The move limit constraint is accounted for through the bounds 
\begin{equation}
\left[\mathbf{z}_{n}^{\mathsf{L}}\right]_{k}=\max\left(\delta_{\rho},\left[\mathbf{z}_{n}\right]_{k}-m_{n}\right),\quad\left[\mathbf{z}_{n}^{\mathsf{U}}\right]_{k}=\min\left(1,\left[\mathbf{z}_{n}\right]_{k}+m_{n}\right),\qquad k=1,\dots,m\label{eq:z_L-z_U}
\end{equation}
In order to embed the curvature information from the reciprocal approximation
(\ref{eq:J_rec_z_y}) in the quadratic model, we choose 
\begin{equation}
\mathbf{H}_{n}=\mbox{diag}(\hat{h}_{1}(\mathbf{z}_{n}),\dots,\hat{h}_{m}(\mathbf{z}_{n}))\label{eq:Hn_z_rec}
\end{equation}
where $\hat{h}_{k}(\mathbf{z}_{n})\equiv\max(h_{k}(\mathbf{z}_{n}),\delta_{E})$
and, as defined before, $0<\delta_{E}\ll\lambda$ is a small positive
constant. This modification not only ensures that $\mathbf{H}_{n}$
is positive definite but also that the eigenvalues of $\mathbf{H}_{n}$
are uniformly bounded above and below, a condition that is useful
for the proof of convergence of the algorithm \citep{Bert-book}.
Observe that for all $\mathbf{z}\in\left[\delta_{\rho},1\right]^{m}$,
\begin{equation}
0\leq h_{k}(\mathbf{z})\leq2\delta_{\rho}^{-1}\left\Vert \mathbf{E}(\mathbf{z})\right\Vert _{\infty}\leq2p\delta_{\rho}^{-p-1}M_{h}\left\Vert \mathbf{U}(\mathbf{z})\right\Vert ^{2}\leq2p\delta_{\rho}^{-p-1}M_{h}c_{h}^{-2}\left\Vert \mathbf{F}\right\Vert ^{2}\label{eq:upp-bnd-h_k}
\end{equation}
where we used the fact that $\mathbf{U}^{T}\mathbf{k}_{e}\mathbf{U}\leq\delta_{\rho}^{-p}\mathbf{U}^{T}\mathbf{K}(\mathbf{z})\mathbf{U}\leq\delta_{\rho}^{-p}M_{h}\left\Vert \mathbf{U}(\mathbf{z})\right\Vert ^{2}$
and that the eigenvalues of $\mathbf{K}^{-1}$ are bounded above by
$c_{h}^{-1}$.

The step size parameter $\tau_{n}$ in (\ref{eq:z-Quad-Subproblem})
must be sufficiently small so that the quadratic model is a conservative
approximation and majorizes $J$. If $\tau_{n}>0$ is chosen so that
the update $\mathbf{z}_{n+1}$ satisfies
\begin{equation}
J(\mathbf{z}_{n+1})\leq Q_{J}(\mathbf{z}_{n+1};\mathbf{z}_{n},\tau_{n})\label{eq:majorization}
\end{equation}
then one can show \citep{Bert-book} 
\begin{equation}
\tilde{J}(\mathbf{z}_{n})-\tilde{J}(\mathbf{z}_{n+1})\geq\frac{1}{2\tau_{n}}\left\Vert \mathbf{z}_{n}-\mathbf{z}_{n+1}\right\Vert _{\mathbf{H}_{n}}^{2}\label{descent}
\end{equation}
If $\mathbf{z}_{n}$ is a stationary point of $\tilde{J}$, that is
$(\mathbf{z}-\mathbf{z}_{n})^{T}\nabla\tilde{J}(\mathbf{z}_{n})\geq0$
for all $\mathbf{z}\in\left[\delta_{\rho},1\right]^{m}$, then $\mathbf{z}_{n+1}=\mathbf{z}_{n}$
for all $\tau_{n}>0$. To see this, we write (\ref{eq:z-Quad-Subproblem})
equivalently as
\begin{equation}
\min_{\mathbf{z}_{n}^{\mathsf{L}}\leq\mathbf{z}\leq\mathbf{z}_{n}^{\mathsf{U}}}(\mathbf{z}-\mathbf{z}_{n})^{T}\nabla\tilde{J}(\mathbf{z}_{n})+\frac{1}{2\tau_{n}}\left\Vert \mathbf{z}-\mathbf{z}_{n}\right\Vert _{\mathbf{H}_{n}+\tau_{n}\mathbf{G}}^{2}\label{eq:proof-zn-stationary}
\end{equation}
Since $\mathbf{H}_{n}+\tau_{n}\mathbf{G}$ is positive definite and
$\mathbf{z}_{n}$ is a stationary point, the objective function is
strictly positive for all $\mathbf{z}\in\left[\mathbf{z}_{n}^{\mathsf{L}},\mathbf{z}_{n}^{\mathsf{U}}\right]$
with $\mathbf{z}\neq\mathbf{z}_{n}$ while it vanishes at $\mathbf{z}_{n}$,
thereby establishing optimality of $\mathbf{z}_{n}$ for subproblem
(\ref{eq:z-Quad-Subproblem}). Otherwise, if $\mathbf{z}_{n}$ is
not a stationary point of $\tilde{J}$, then\emph{ $\mathbf{z}_{n+1}\neq\mathbf{z}_{n}$}
for sufficiently small $\tau_{n}$, and (\ref{descent}) shows that
there is a decrease in the objective function. This latter fact shows\emph{
}that\emph{ the algorithm is monotonically decreasing}.

A step size parameter satisfying (\ref{eq:majorization}) is guaranteed
to exist if $J$ has a Lipschitz gradient, that is, for some positive
constant $L$, 
\begin{equation}
\left\Vert \nabla J(\mathbf{z})-\nabla J(\mathbf{y})\right\Vert \leq L\left\Vert \mathbf{z}-\mathbf{y}\right\Vert ,\quad\forall\mathbf{z},\mathbf{y}\in\mathrm{dom}(J)\label{eq:grad_J_Lipschitz}
\end{equation}
One can show%
\footnote{This is in fact stronger than (\ref{eq:majorization})%
} $J(\mathbf{z})\leq Q_{J}(\mathbf{z};\mathbf{z}_{n};\tau_{n})$ for
all $\mathbf{z}\in\left[\delta_{\rho},1\right]^{m}$ if the step size
satisfies 
\begin{equation}
\tau_{n}^{-1}\mathbf{H}_{n}>L\mathbf{I}\label{eq:tau_n_Lipschitz}
\end{equation}
in the sense of quadratic forms, i.e., $\tau_{n}^{-1}\mathbf{H}_{n}-L\mathbf{I}$
is positive definite \citep{Bert-book}. We verify that the gradient
of compliance $\nabla J$ given by (\ref{eq:Grad_J}) is indeed Lipschitz:
\begin{eqnarray}
\left\Vert \nabla J(\mathbf{z})-\nabla J(\mathbf{y})\right\Vert  & = & p\left\Vert \mathbf{E}(\mathbf{z})-\mathbf{E}(\mathbf{y})\right\Vert \nonumber \\
 & \leq & p\left[\sum_{e=1}^{l}\left(\left[\mathbf{P}\mathbf{z}\right]_{e}^{p-1}\delta_{\rho}^{-p}M_{h}\left\Vert \mathbf{U}(\mathbf{z})\right\Vert ^{2}-\left[\mathbf{P}\mathbf{y}\right]_{e}^{p-1}\delta_{\rho}^{-p}M_{h}\left\Vert \mathbf{U}(\mathbf{y})\right\Vert ^{2}\right)^{2}\right]^{1/2}\nonumber \\
 & \leq & p\delta_{\rho}^{-p}M_{h}\left[\sum_{e=1}^{l}\left(\left[\mathbf{P}\mathbf{z}\right]_{e}c_{h}^{-2}\left\Vert \mathbf{F}\right\Vert ^{2}-\left[\mathbf{P}\mathbf{y}\right]_{e}c_{h}^{-2}\left\Vert \mathbf{F}\right\Vert ^{2}\right)^{2}\right]^{1/2}\label{eq:nabla_J_Lips_proof}\\
 & \leq & p\delta_{\rho}^{-p}M_{h}c_{h}^{-2}\left\Vert \mathbf{F}\right\Vert ^{2}\left\Vert \mathbf{P}\mathbf{z}-\mathbf{P}\mathbf{y}\right\Vert \nonumber \\
 & \leq & p\delta_{\rho}^{-p}M_{h}c_{h}^{-2}\left\Vert \mathbf{F}\right\Vert ^{2}\left\Vert \mathbf{z}-\mathbf{y}\right\Vert \nonumber 
\end{eqnarray}
The step size $\tau_{n}$ can be selected with \emph{a priori} knowledge
of the Lipschitz constant $L$ but this may be too conservative and
may slow down the convergence of the algorithm. Instead, in each iteration,
one can gradually decrease the step size via a backtracking routine
until $\mathbf{z}_{n+1}$ satisfies (\ref{eq:majorization}). An alternative,
possibly weaker, descent condition is the Armijo rule which requires
that for some constant $0<\nu<1$, the update satisfies
\begin{equation}
\tilde{J}(\mathbf{z}_{n})-\tilde{J}(\mathbf{z}_{n+1})\geq\nu\left(\mathbf{z}_{n}-\mathbf{z}_{n+1}\right)^{T}\nabla\tilde{J}(\mathbf{z}_{n})\label{eq:Armijo}
\end{equation}
Though the implementation of such step size routines is straightforward,
due to the high cost of function evaluations for the compliance problem
(which requires solving the state equation to compute the value of
$J$), the number of trials in satisfying the descent condition must
be limited. Therefore, there is a tradeoff between attempting to choose
a large step size to speed up convergence and the cost associated
with the selection routine. As shown in the next section, we have
found that fixing $\tau_{n}=1$, which eliminates the cost of backtracking
routine, generally leads to a stable and convergent algorithm. In
some cases, however, the overall cost can be reduced by using larger
step sizes.

As in section 3, ignoring constant terms in $\mathbf{z}_{n}$ and
rearranging, we can write (\ref{eq:z-Quad-Subproblem}) equivalently
as
\begin{equation}
\mathbf{z}_{n+1}=\underset{\mathbf{z}_{n}^{\mathsf{L}}\leq\mathbf{z}\leq\mathbf{z}_{n}^{\mathsf{U}}}{\mbox{argmin}}\ \left\Vert \mathbf{z}-\mathbf{z}_{n+1}^{*}\right\Vert _{\mathbf{H}_{n}+\tau_{n}\mathbf{G}}^{2}\label{eq:disc_FBS_alt_form}
\end{equation}
where the interim update $\mathbf{z}_{n+1}^{*}$ is the given by
\begin{equation}
\mathbf{z}_{n+1}^{*}=\mathbf{z}_{n}-\tau_{n}\left(\mathbf{H}_{n}+\tau_{n}\mathbf{G}\right)^{-1}\left[\nabla\tilde{J}(\mathbf{z}_{n})\right]\label{eq:z_n+1_star}
\end{equation}
With the appropriate choice of step size (satisfying any one of the
conditions (\ref{eq:majorization}), (\ref{eq:tau_n_Lipschitz}),
or (\ref{eq:Armijo})) and boundedness of $\mathbf{H}_{n}$, it can
be shown that every limit point of the the sequence $\mathbf{z}_{n}$
generated by the algorithm is a critical point of $\tilde{J}$. For
the particular case of quadratic regularization, it is evident from
(\ref{eq:z_n+1_star}) that the algorithm reduces to the so-called
scaled gradient projection algorithm, and the convergence proof can
be found in \citep{Bert-book}. A more general proof can be found
in the review paper on proximal splitting method by \citep{Beck-review}
though the metric associated with the proximal term, i.e., $\left\Vert \mathbf{z}-\mathbf{z}_{n}\right\Vert _{\mathbf{H}_{n}+\tau_{n}\mathbf{G}}^{2}$
in (\ref{eq:z-Quad-Subproblem}), is fixed there.

As seen from (\ref{eq:z-Quad-Subproblem}) or (\ref{eq:disc_FBS_alt_form}),
the forward-backward algorithm requires the solution to a sparse,
strictly convex quadratic program subject to simple bound constraints
which can be efficiently solved using a variety of methods, e.g.,
the active set method. Alternatively, the projection of $\mathbf{z}_{n+1}^{*}$
can be recast as a bound constrained sparse least squares problem
and solved using algorithms in \citep{Adlers-thesis}.

\subsection*{Two-metric projection variation}

Next we discuss a variation of the splitting algorithm that simplifies
the projection step (\ref{eq:disc_FBS_alt_form}) by augmenting the
interim density (\ref{eq:z_n+1_star}). More specifically, we adopt
a variant of the two-metric projection method \citep{Bert-ProjNewton,Gafni-TMP},
in which the norm in (\ref{eq:disc_FBS_alt_form}) is replaced by
the usual Euclidean norm, and the scaling matrix $\mathbf{H}_{n}+\tau_{n}\mathbf{G}$
in the interim step (\ref{eq:z_n+1_star}) is made diagonal with respect
to the active components of $\mathbf{z}_{n}$. 

Let $I_{n}=I_{n}^{\mathsf{L}}\cup I_{n}^{\mathsf{U}}$ denote the
set of active constraints where
\begin{eqnarray}
I_{n}^{\mathsf{L}} & = & \left\{ k:\left[\mathbf{z}_{n}\right]_{k}\leq\delta_{\rho}+\epsilon\mbox{ and }\left[\nabla\tilde{J}(\mathbf{z}_{n})\right]_{k}>0\right\} \label{eq:I_n^L}\\
I_{n}^{\mathsf{U}} & = & \left\{ k:\left[\mathbf{z}_{n}\right]_{k}\geq1-\epsilon\mbox{ and }\left[\nabla\tilde{J}(\mathbf{z}_{n})\right]_{k}<0\right\} \label{eq:I_n^U}
\end{eqnarray}
Here $\epsilon$ is an algorithmic parameter (we fix it at $10^{-3}$
for the numerical results) that enlarges the set of active constraints
in order to avoid the discontinuities that may otherwise arise \citep{Bert-ProjNewton}.
Then 
\begin{equation}
\left[\mathbf{D}_{n}\right]_{ij}\equiv\begin{cases}
0 & \mbox{if }i\neq j\mbox{ and }i\in I_{n}\mbox{ or }j\in I_{n}\\
\left[\mathbf{H}_{n}+\tau_{n}\mathbf{G}\right]_{ij} & \mbox{otherwise}
\end{cases}\label{eq:D_n}
\end{equation}
is a scaling matrix formed from $\mathbf{H}_{n}+\tau_{n}\mathbf{G}$
that is diagonal with respect to $I_{n}$ and therefore removes the
coupling between the active and free constraints. The operation in
(\ref{eq:D_n}) essentially consists of zeroing out all the off-diagonal
entries of $\mathbf{H}_{n}+\tau_{n}\mathbf{G}$ for the active components.
Note that any other positive matrix with the same structure as $\mathbf{D}_{n}$
can be used. The new interim density is then defined as 
\begin{equation}
\mathbf{z}_{n+1}^{*}=\mathbf{z}_{n}-\tau_{n}\mathbf{D}_{n}^{-1}\left[\nabla\tilde{J}(\mathbf{z}_{n})\right]\label{eq:TMP_zn+1_star}
\end{equation}
and the next iterate is given by the Euclidian projection of this
interim density onto the constraint set 
\begin{equation}
\mathbf{z}_{n+1}=\underset{\mathbf{z}_{n}^{\mathsf{L}}\leq\mathbf{z}\leq\mathbf{z}_{n}^{\mathsf{U}}}{\mbox{argmin}}\ \left\Vert \mathbf{z}-\mathbf{z}_{n+1}^{*}\right\Vert ^{2}\label{eq:TMP_L2_proj}
\end{equation}
which has an explicit solution 
\begin{equation}
\left[\mathbf{z}_{n+1}\right]_{k}=\min\left(\max\left(\left[\mathbf{z}_{n}^{\mathsf{L}}\right]_{k},\left[\mathbf{z}_{n+1}^{*}\right]_{k}\right),\left[\mathbf{z}_{n}^{\mathsf{U}}\right]_{k}\right),\quad k=1,\dots,m\label{eq:TMP_zn+1}
\end{equation}
Since $\mathbf{D}_{n}^{-1}\nabla\tilde{J}(\mathbf{z}_{n})$ can be
viewed as the gradient of $\tilde{J}$ with respect to the metric
induced by $\mathbf{D}_{n}$, we can see that the present algorithm
consisting of (\ref{eq:TMP_zn+1_star}) and (\ref{eq:TMP_L2_proj})
utilizes two separate metrics for differentiation and projection operations.
The significant computational advantage of carrying out the projection
step with respect to the Euclidian norm is due to the particular separable
structure of the constraint set. Compared to the forward-backward
algorithm discussed before, at the cost of modifying the scaling matrix,
the overhead associated with solving the quadratic program (cf. (\ref{eq:disc_FBS_alt_form}))
is eliminated.
\begin{figure}
\centering{}\includegraphics[scale=0.6]{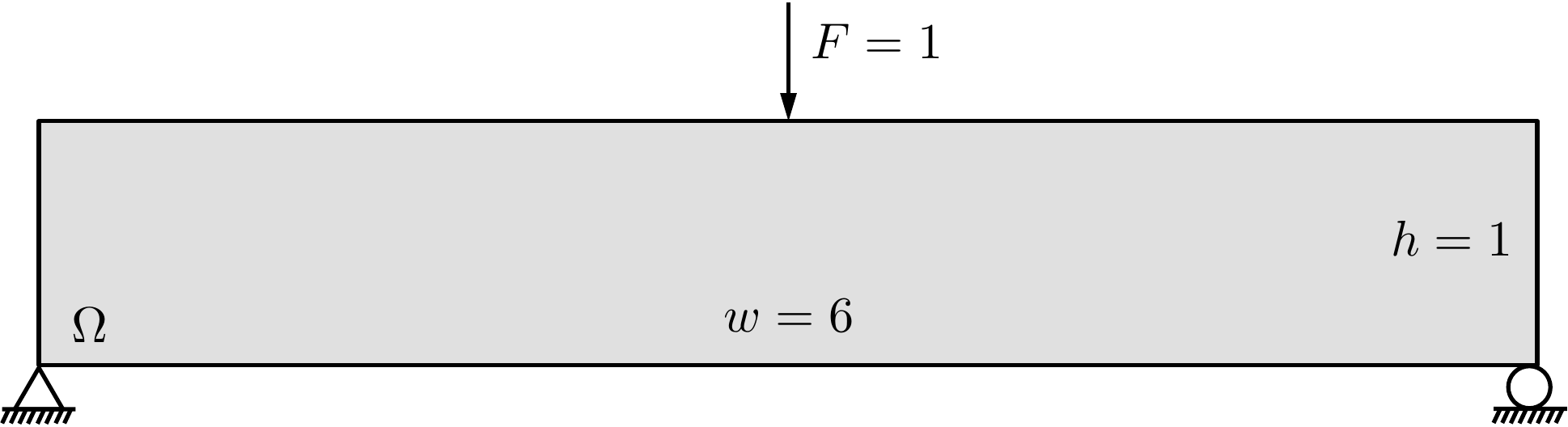}\caption{The design domain and boundary conditions for the MBB beam problem
\label{fig:MBB-domain}}
\end{figure}

As in the previous algorithm, one can show that $\mathbf{z}_{n}$
is a critical point of $\tilde{J}$ if and only if $\mathbf{z}_{n+1}=\mathbf{z}_{n}$
for all $\tau_{n}>0$. Similarly, if $\mathbf{z}_{n}$ is not a stationary
point, then for a sufficiently small step size, the next iterate decreases
the value of the cost function, i.e., $\tilde{J}(\mathbf{z}_{n+1})<\tilde{J}(\mathbf{z}_{n})$.
The choice of $\tau_{n}$ can be again obtained from an Amijo-type
condition along the projection arc (cf. \citep{Bert-ProjNewton}),
namely,
\begin{equation}
\tilde{J}(\mathbf{z}_{n})-\tilde{J}(\mathbf{z}_{n+1})\geq\nu\mathbf{d}_{n}^{T}\nabla\tilde{J}(\mathbf{z}_{n})\label{eq:Armijo_TMP}
\end{equation}
where the direction vector $\mathbf{d}_{n}$ is given by 
\begin{equation}
\left[\mathbf{d}_{n}\right]_{k}=\begin{cases}
\left[\mathbf{z}_{n}\right]_{k}-\left[\mathbf{z}_{n+1}\right]_{k} & k\in I_{n}\\
\left[\tau_{n}\mathbf{D}_{n}^{-1}\nabla\tilde{J}(\mathbf{z}_{n})\right]_{k} & k\notin I_{n}
\end{cases}\label{eq:d_n_for_TMP}
\end{equation}
In the next section, we will compare the performance of the forward-backward
algorithm consisting of (\ref{eq:disc_FBS_alt_form}) and (\ref{eq:z_n+1_star})
with the two-metric projection consisting of (\ref{eq:TMP_zn+1_star})
and (\ref{eq:TMP_zn+1}).

\section{Numerical Investigations}

The model compliance minimization problem adopted here is the benchmark
MBB beam problem, whose domain geometry and prescribed loading and
boundary conditions are shown in Figure \ref{fig:MBB-domain}. Using
appropriate boundary conditions, the symmetry of the problem is exploited
to pose and solve the state equation only on half of the extended
domain. The constituent material $\mathbf{C}_{0}$ is assumed to be
isotropic with unit Young's modulus and Poisson ratio of $0.3$. The
volume penalty parameter is $\lambda=200/\left|\Omega\right|$ where
$\left|\Omega\right|$ is the area of the extended design domain.
For all the results in this section, the lower bound on the density
is set to $\delta_{\rho}=10^{-3}$ and, unless otherwise stated, the
SIMP penalty exponent is fixed at $p=3$. A simple backtracking algorithm
is used to determine the value of the step size parameter. Given constants
$\tau_{0}>0$ and $0<\sigma<1$, the step size parameter in the $n$th
iteration is given by{\footnotesize{}}
\begin{table}
\begin{centering}
{\footnotesize{}}%
\begin{tabular}{|c|c|c||c|c|c|c|c|c|c|c|}
\hline 
{\footnotesize{algorithm}} & {\footnotesize{$\mathbf{H}_{n}$}} & {\footnotesize{$\tau_{0}$}} & {\footnotesize{\# it.}} & {\footnotesize{\# bt.}} & {\footnotesize{$\ell(\mathbf{u}_{\rho})$}} & {\footnotesize{$R(\rho)$}} & {\footnotesize{$V(\rho)$}} & {\footnotesize{$\tilde{J}(\rho)$}} & {\footnotesize{$E_{1}$}} & {\footnotesize{$E_{2}$}}\tabularnewline
\hline 
\hline 
{\footnotesize{FBS}} & {\footnotesize{identity}} & {\footnotesize{1}} & {\footnotesize{316}} & {\footnotesize{0}} & {\footnotesize{100.019}} & {\footnotesize{8.553}} & {\footnotesize{0.5120}} & {\footnotesize{210.965}} & {\footnotesize{9.962e-6}} & {\footnotesize{9.943e-5}}\tabularnewline
\hline 
{\footnotesize{FBS}} & {\footnotesize{identity}} & {\footnotesize{2}} & {\footnotesize{215}} & {\footnotesize{154}} & {\footnotesize{100.093}} & {\footnotesize{8.537}} & {\footnotesize{0.5114}} & {\footnotesize{210.914}} & {\footnotesize{9.178e-6}} & {\footnotesize{5.812e-5}}\tabularnewline
\hline 
{\footnotesize{FBS}} & {\footnotesize{reciprocal}} & {\footnotesize{1}} & {\footnotesize{186}} & {\footnotesize{0}} & {\footnotesize{99.937}} & {\footnotesize{8.594}} & {\footnotesize{0.5125}} & {\footnotesize{211.032}} & {\footnotesize{9.769e-6}} & {\footnotesize{9.363e-5}}\tabularnewline
\hline 
{\footnotesize{FBS}} & {\footnotesize{reciprocal}} & {\footnotesize{2}} & {\footnotesize{91}} & {\footnotesize{39}} & {\footnotesize{100.095}} & {\footnotesize{8.568}} & {\footnotesize{0.5117}} & {\footnotesize{211.008}} & {\footnotesize{4.926e-6}} & {\footnotesize{9.746e-5}}\tabularnewline
\hline 
\hline 
{\footnotesize{TMP}} & {\footnotesize{identity}} & {\footnotesize{1}} & {\footnotesize{330}} & {\footnotesize{0}} & {\footnotesize{100.076}} & {\footnotesize{8.533}} & {\footnotesize{0.5117}} & {\footnotesize{210.951}} & {\footnotesize{9.958e-6}} & {\footnotesize{9.973e-5}}\tabularnewline
\hline 
{\footnotesize{TMP}} & {\footnotesize{identity}} & {\footnotesize{2}} & {\footnotesize{151}} & {\footnotesize{78}} & {\footnotesize{100.060}} & {\footnotesize{8.556}} & {\footnotesize{0.5116}} & {\footnotesize{210.938}} & {\footnotesize{9.639e-6}} & {\footnotesize{5.900e-5}}\tabularnewline
\hline 
{\footnotesize{TMP}} & {\footnotesize{reciprocal}} & {\footnotesize{1}} & {\footnotesize{179}} & {\footnotesize{0}} & {\footnotesize{99.943}} & {\footnotesize{8.592}} & {\footnotesize{0.5125}} & {\footnotesize{211.031}} & {\footnotesize{9.878e-6}} & {\footnotesize{9.453e-5}}\tabularnewline
\hline 
{\footnotesize{TMP}} & {\footnotesize{reciprocal}} & {\footnotesize{2}} & {\footnotesize{85}} & {\footnotesize{34}} & {\footnotesize{100.078}} & {\footnotesize{8.578}} & {\footnotesize{0.5117}} & {\footnotesize{210.999}} & {\footnotesize{9.043e-6}} & {\footnotesize{8.074e-5}}\tabularnewline
\hline 
\end{tabular}
\par\end{centering}{\footnotesize \par}

{\footnotesize{\caption{{\footnotesize{Summary of influence of various factors in the algorithm
for the MBB problem with $\beta=0.06$. The acronym FBS designates
the forward-backward algorithm and TMP refers to the two-metric projection
algorithm. Forth and fifth columns show the total number of iterations
and backtracking steps. The remaining columns show the final value
of compliance $\ell(\mathbf{u}_{\rho})$, regularization term $R(\rho)$,
volume fraction $V(\rho)=\left|\Omega\right|^{-1}\int_{\Omega}\rho\mathrm{d}\mathbf{x}$,
the regularized objective $\tilde{J}(\rho)$, the relative change
in cost function value $E_{1}$ and the error in satisfaction of the
first order conditions of optimality $E_{2}$ \label{tab:Summary-Mbb-0-06}}}}
}}
\end{table}
\begin{equation}
\tau_{n}=\sigma^{k_{n}}\tau_{0}\label{eq:Backtracking}
\end{equation}
where $k_{n}$ is the smallest non-negative integer such that $\tau_{n}$
satisfies (\ref{eq:Armijo}) or (\ref{eq:Armijo_TMP}). In practice,
this means that we begin with the initial step size $\tau_{0}$ and
reduce it by a factor of $\sigma$ until descent conditions are satisfied.
The descent parameter is set to $\nu=10^{-3}$ and the backtracking
parameter is $\sigma=0.6$. Note that larger $\nu$ leads to a more
severe descent requirement and subsequently smaller $\tau_{n}$. Similarly,
smaller $\sigma$ reduces the step size parameter by a larger factor
which can decrease the number of backtracking step. Note, however,
that using small step sizes may lead to slow convergence of the algorithm.

Since each backtracking step involves evaluating the cost functional
and therefore solving the state equation, as a measure of computational
cost, we keep track of the total number of backtracking steps (i.e.,
$\sum_{n}k_{n}$) in addition to the total number of iterations. The
convergence criteria adopted here is based on the relative decrease
in the objective function
\begin{equation}
E_{1}=\frac{\bigl|\tilde{J}(\mathbf{z}_{n+1})-\tilde{J}(\mathbf{z}_{n})\bigr|}{\bigl|\tilde{J}(\mathbf{z}_{n})\bigr|}\leq\epsilon_{1}\label{eq:Error1}
\end{equation}
and the satisfaction of the first order conditions of optimality according
to 
\begin{equation}
E_{2}=\frac{\bigl\Vert\mathcal{P}[\mathbf{z}_{n+1}-\nabla\tilde{J}(\mathbf{z}_{n+1})]-\mathbf{z}_{n+1}\bigr\Vert}{\bigl\Vert\mathbf{z}_{n+1}\bigr\Vert}\leq\epsilon_{2}\label{eq:Error2}
\end{equation}
Here $\mathcal{P}$ is the Euclidian projection onto the constraint
set $\left[\delta_{\rho},1\right]^{m}$ defined by $\left[\mathcal{P}(\mathbf{y})\right]_{i}=\min\left(\max\left(0,\left[\mathbf{y}\right]_{i}\right),1\right)$.
Unless otherwise stated, we have selected $\epsilon_{1}=10^{-5}$
and $\epsilon_{2}=10^{-4}$.

We begin with the investigation of the behavior of two forms of the
algorithm with different choice of parameters discussed in the previous
section. In particular, we compare the forward-backward algorithm
with the two-metric projection method and investigate the influence
of the Hessian approximation. In addition to the choice of $\mathbf{H}_{n}$
defined by (\ref{eq:Hn_z_rec}), we also consider a fixed scaling
of the identity matrix 
\begin{figure}
\begin{centering}
\includegraphics[scale=0.6]{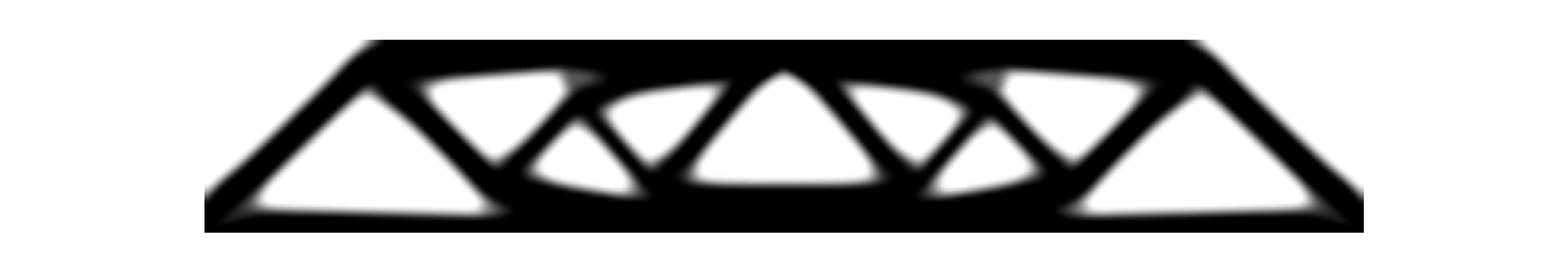}
\par\end{centering}

\caption{Final density field for the MBB problem and $\beta=0.06$ plotted
in grayscale. This result was generated using the TMP algorithm with
$\tau_{0}=2$ and $m_{n}=1$ \label{fig:Mbb-result-06}}
\end{figure}
 
\begin{equation}
\mathbf{H}_{n}\equiv\alpha\mathbf{I},\qquad n=1,2,\dots\label{eq:constant-Hessian}
\end{equation}
for which the algorithm becomes the basic forward-backward algorithm
with the same proximal term in every iteration. The scaling coefficient
$\alpha$ is set to $4\lambda A$ where $A$ is the area of an element.
This choice is made so that the step size parameter $\tau_{n}$ is
the same order of magnitude as with reciprocal Hessian. The other
parameter investigated here is the initial step size parameter $\tau_{0}$
and we consider two choices $\tau_{0}=1$ and $\tau_{0}=2$. In all
cases, the move limit is fixed at $m_{n}=1$ for all $n$ and thus
$\mathcal{A}_{n}=\mathcal{A}$. 

The model problem is the MBB beam discretized with a grid of 300 by
50 bilinear quad elements and Tikhonov regularization parameter is
set to $\beta=0.06$. The initial guess in all cases is taken to be
uniform density field $\rho_{h}\equiv1/2$. All the possible combinations
of the above choices produce the same final topology, similar to the
representative solution shown in Figure \ref{fig:Mbb-result-06}.
\emph{This shows the framework exhibits stable convergence to the
same final solution and is relatively insensitive to various choices
of algorithmic parameters for this level of regularization.} What
is different, however, is the speed of convergence and the required
computational effort as measured by the number of the backtracking
steps, total number of iterations, and cost per iteration. The results
are summarized in Table \ref{tab:Summary-Mbb-0-06}.

First we note that the initial step size $\tau_{0}=1$ does not lead
to any backtracking steps which means that in each iteration the step
size parameter is $\tau_{n}=1$. By contrast, using the larger initial
step size parameter $\tau_{0}=2$ requires backtracking steps to satisfy
the descent condition but substantially reduces the total number of
iterations. Moreover, in all cases, the constant Hessian (\ref{eq:constant-Hessian})
requires nearly twice as many iterations and backtracking steps compared
to the ``reciprocal'' Hessian. \emph{This highlights the fact that
embedding the reciprocal approximation of compliance does indeed lead
to faster convergence.} Overall, the best performance is obtained
using the reciprocal approximation and larger initial step size parameter. 

For this problem, the forward-backward algorithm and the two-metric
projection method roughly have the same number of iterations and backtracking
steps. However, the cost per iteration for the two-metric projection
is significantly lower since the projection step is computationally
trivial. Therefore, the two-metric projection is more efficient. 
\begin{figure}
\begin{centering}
\includegraphics[scale=0.6]{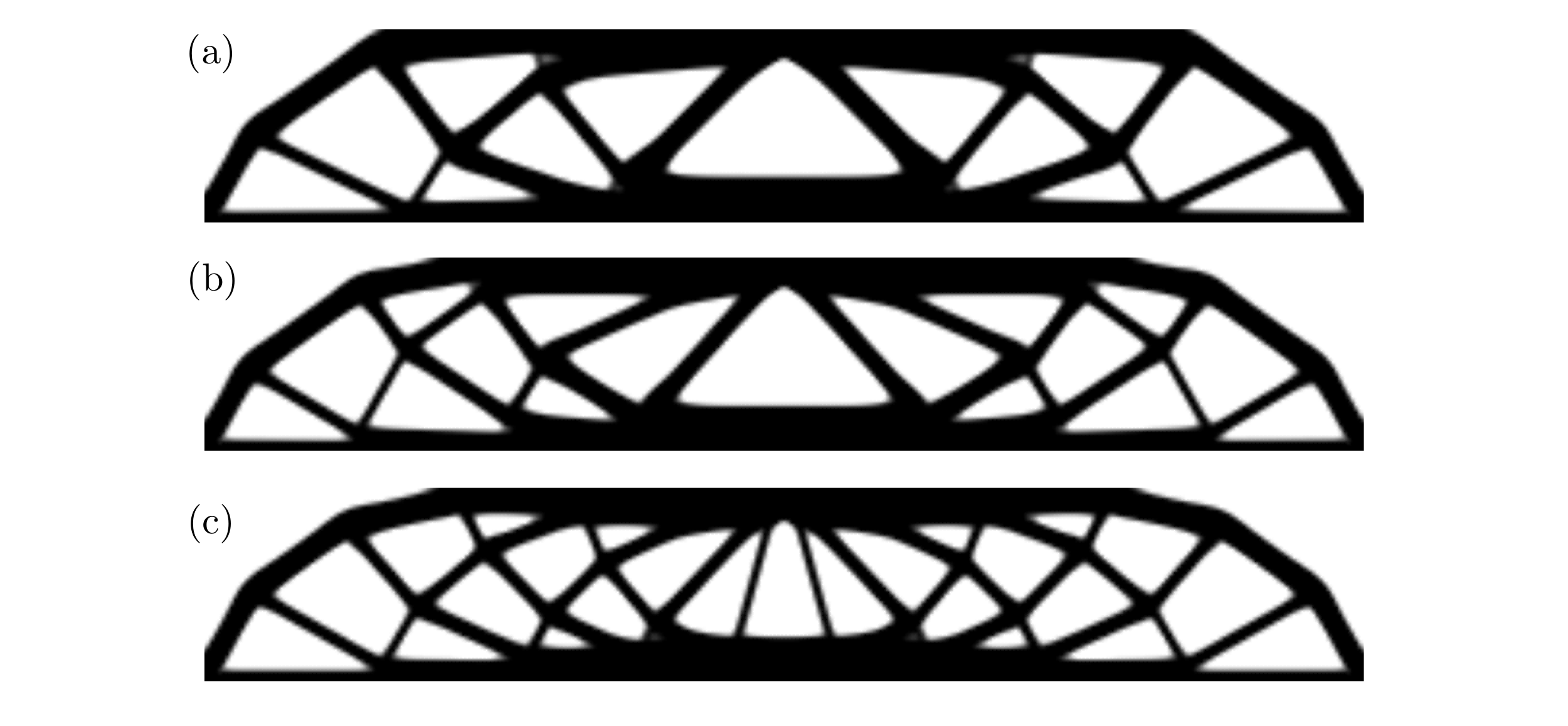}
\par\end{centering}

\caption{Final densities plotted in grayscale for the MBB problem and $\beta=0.01$.
The results are generated using the TMP algorithm with (a) $\tau_{0}=2$,
$m_{n}=1$ (b) $\tau_{0}=1$, $m_{n}=1$ (c) $\tau_{0}=1$, $m_{n}=0.03$
\label{fig:Mbb-01}}
\end{figure}

Next we investigate the performance of the algorithm for a smaller
value of the regularization parameter which is expected to produce
more complex topologies. For the next set of results, we set $\beta=0.01$.
In all cases considered, the forward-backward and the two-metric projection
algorithms both give identical final topologies with roughly the same
number of iterations and so we only report the results for the two-metric
projection algorithm. Also, as demonstrated by the first study, the
use of reciprocal approximation leads to better and faster convergence
of the algorithm so we limit the remaining results to the ``reciprocal''
$\mathbf{H}_{n}$. The tolerance level $\epsilon_{2}=10^{-4}$ for
satisfaction of the optimality condition is relatively stringent in
this case due to the complexity of final designs (compared to $\beta=0.06$)
and leads to a large number of iterations with little change in density
near the optimum. We therefore increase the tolerance to $\epsilon_{2}=2\times10^{-4}$
which gives nearly identical final topologies but with fewer iterations.

We examine the influence of the step size parameter and move limit,
which unlike the previous case of large regularization parameter,
can lead to different final solutions. We consider two possible initial
step size parameters $\tau_{0}=1$ and $\tau_{0}=2$, as well as two
choices for the move limit $m_{n}\equiv1$ and $m_{n}\equiv0.03$.
Here we are using a fixed move limit $m_{n}$ for all iteration $n$.
It may be possible to devise a strategy to increase $m_{n}$ in the
later stages of optimization to improve convergence. The results are
summarized in Table \ref{tab:Summary-Mbb-01} and the final solutions
are shown in Figure \ref{fig:Mbb-01}. 

First note that with no move limit constraints, i.e., $m_{n}=1$,
the final solution with the more aggressive choice of initial step
size parameter ($\tau_{0}=2$) is less complex and has fewer members
compared to $\tau_{0}=1$, which as before does not require any backtracking
steps. Note, however, that the more aggressive scheme in fact requires
more iterations to converge. In the presence of move limits, there
is no backtracking step with either choice of step size but the larger
step size does reduce the total number of iterations. The final topologies
are identical and have more members compared to the solutions obtained
without the move limits. It is interesting to note that the overall
iteration count is lowest for $\tau_{0}=2$ and $m_{n}=0.03$ despite
the limit on the change in density in each iteration. As noted earlier,
the use of move limits can stabilize the convergence of the topology
optimization problem. 
\begin{figure}
\centering{}\includegraphics[scale=0.6]{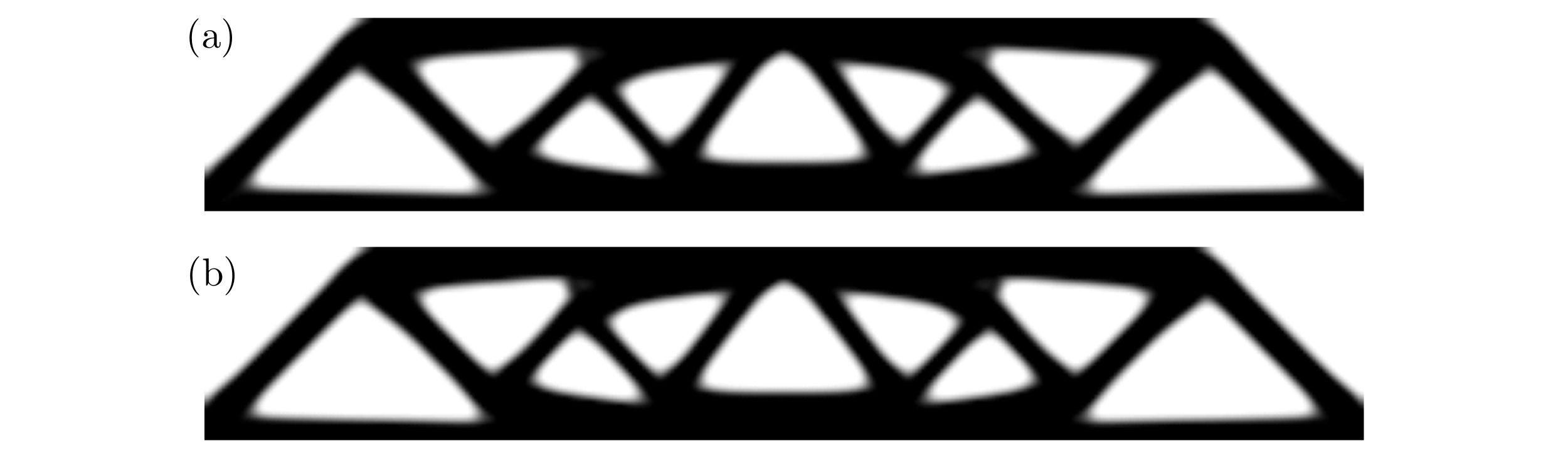}\caption{Final densities plotted in grayscale for the MBB problem with $\beta=0.06$
and SIMP penalty exponent (a) $p=4$ (b) $p=5$\label{fig:HigherSimpP}}
\end{figure}

The overall trend that the move aggressive choice of parameters produce
less complex final solutions is due to the fact that member formation
occurs early on in the algorithm. The most aggressive algorithm ($\tau_{0}=2$,
$m_{n}=1$) still produces the best solution as measured by $\tilde{J}$
while the solution obtained enforcing the move limit $m_{n}=0.03$
has the lowest value of compliance $J$ (due to distribution of members
and slightly higher volume fraction).

We note that aside from the higher degree of complexity, the optimal
densities for $\beta=0.01$ contain fewer intermediate values compared
to the solution for $\beta=0.06$. One measure of discreteness used
in \citep{Sigmund:2007p460} is given by
\begin{equation}
M(\rho)=\frac{1}{\left|\Omega\right|}\int_{\Omega}4\left(\rho-\delta_{\rho}\right)\left(1-\rho\right)\mathrm{d}\mathbf{x}
\end{equation}
which is equal to zero if $\rho$ takes only values of $\delta_{\rho}$
and 1. For the solutions shown in Figure \ref{fig:Mbb-01}, $M(\rho)$
is equal to 6.98\%, 7.64\% and 8.90\% from top to bottom, respectively.
In contrast, the optimal density for $\beta=0.06$ (cf. Figure \ref{fig:Mbb-result-06})
has a discreteness measure of 15.0\%. By increasing the value of the
SIMP exponent $p$, the optimal densities can be made more discrete.
The results for $\beta=0.06$ using $p=4$ and $p=5$ are shown in
Figure \ref{fig:HigherSimpP}. While the optimal topologies are nearly
identical to that the solution for $p=3$, the discrete measure is
lower to 13.1\% and 12.1\%, respectively. Observe, however, that the
layer of intermediate densities around the boundary cannot be completely
eliminated even when $p$ is increase to a very large value since
the Tikhonov regularizer is singular in the discontinuous limit of
density. 

\begin{figure}
\begin{centering}
\includegraphics[scale=0.65]{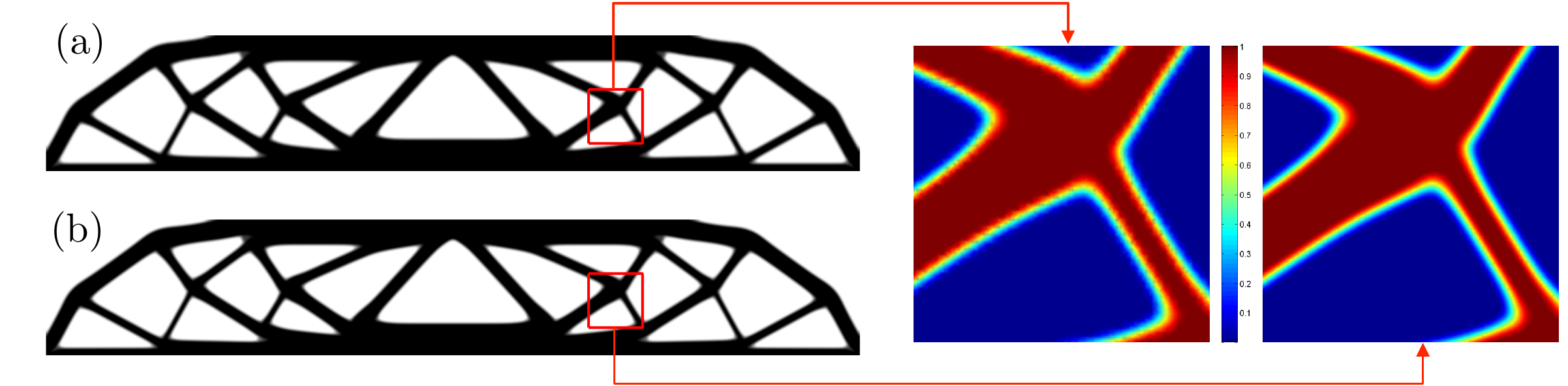}
\par\end{centering}

\caption{Results of the mesh refinement study with (a) $600\times100$ (b)
$1200\times200$ elements. \label{fig:MeshRef} }
\end{figure}

As shown in the previous section, the optimal solutions to the discrete
problem converge to an optimal solution of the continuum problem as
the finite element mesh is refined. We next demonstrate numerically
that solutions produced by the present optimization algorithms appear
to be stable with respect to mesh refinement. We do this for the case
of $\beta=0.01$ using the two-metric projection algorithm with $\tau_{n}\equiv1$
where the final topology is relatively complex and the algorithm is
expected to be more sensitive. As shown in Figure \ref{fig:MeshRef},
we solve the problem using finer grids consisting of $600\times100$
and $1200\times200$ bilinear square elements. The final density distribution
is nearly identical indicating convergence of optimal densities in
the $L^{p}$-norm.

\subsection*{Compliant mechanism design}

The discussion so far has been limited to the problem of compliance
minimization which, as noted earlier, is self-adjoint and its gradient
has the same sign. We conclude this section with design of a compliant
force inverter for which the cost functional is no longer self-adjoint
and therefore, unlike compliance, the gradient field may take both
negative and positive values in the domain. 

\begin{figure}
\begin{centering}
\includegraphics[scale=0.55]{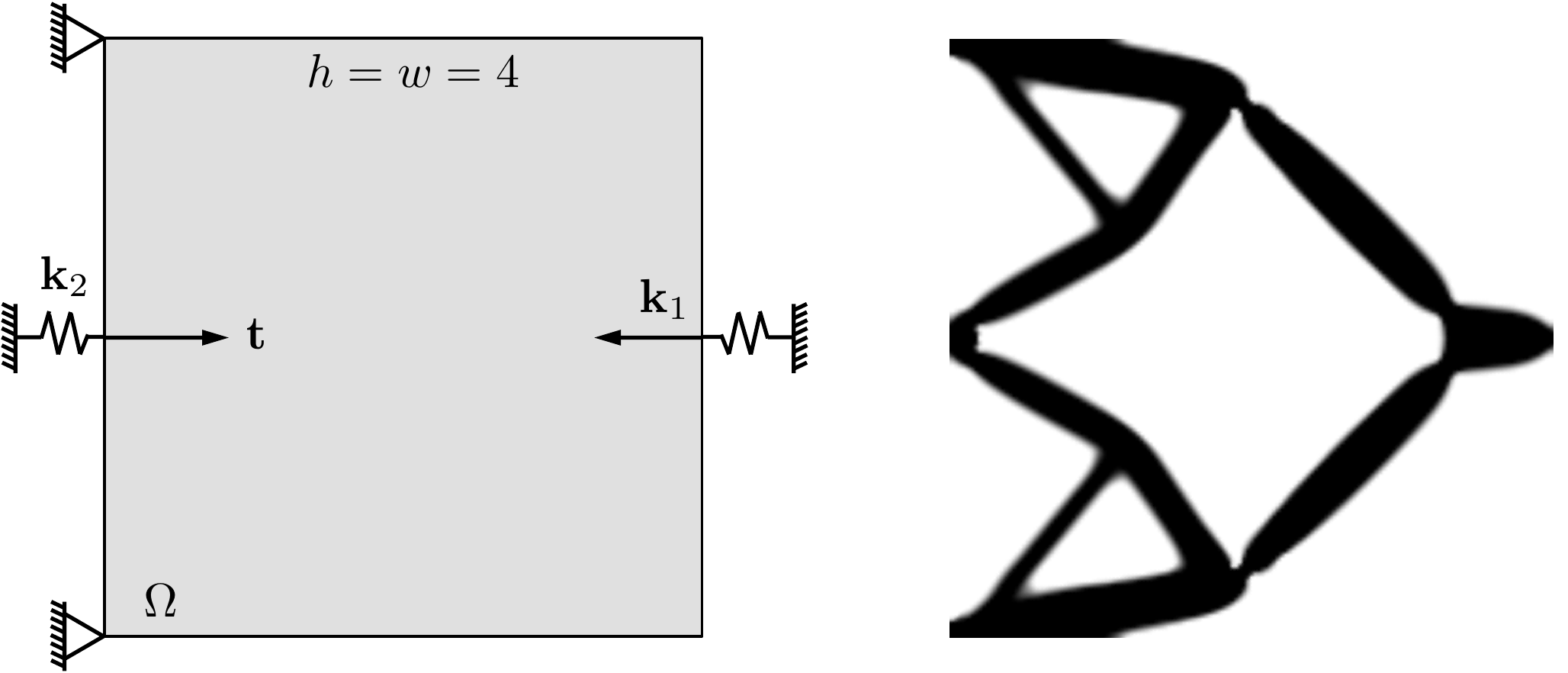}
\par\end{centering}

\caption{The design domain and boundary conditions for the force inverter problem
(left) and the optimal topology (right). For this example, $\left\Vert \mathbf{k}_{1}\right\Vert =\left\Vert \mathbf{k}_{2}\right\Vert =0.1$
\label{fig:Mechanism}}
\end{figure}

The objective of mechanism design is to identify a structure that
maximizes the force exerted on a workpiece under the action of an
external actuator. As illustrated in Figure \ref{fig:Mechanism},
the force inverter transfers the input force of the actuator to a
force at the prescribed output location in the opposite direction.
We assume in this setting that both the workpiece and the actuator
are elastic and their stiffness are represented by vector fields $\mathbf{k}_{1}\in L^{\infty}(\Gamma_{S_{1}})$
and $\mathbf{k}_{2}\in L^{\infty}(\Gamma_{S_{2}})$, respectively.
Here $\Gamma_{S_{1}},\Gamma_{S_{2}}$ are segments of the traction
boundary $\Gamma_{N}\subseteq\partial\Omega$ where the structure
is interacting with these elastic bodies. The tractions experienced
by the structure through this interaction for a displacement field
$\mathbf{u}$ can be written as
\begin{equation}
\mathbf{t}_{S_{r}}(\mathbf{u})=-\left(\mathbf{k}_{r}\cdot\mathbf{u}\right)\frac{\mathbf{k}_{r}}{\left\Vert \mathbf{k}_{r}\right\Vert },\qquad\mbox{on }\Gamma_{S_{r}}\mbox{ for }r=1,2\label{eq:t_k_force}
\end{equation}
Accordingly, the displacement $\mathbf{u}_{\rho}$ for a given distribution
of material $\rho$ in $\Omega$ is the solution to the following
boundary problem
\begin{equation}
a(\mathbf{u}_{\rho},\mathbf{v};\rho)+a_{s}(\mathbf{u}_{\rho},\mathbf{v})=\ell(\mathbf{v})\quad\forall\mathbf{v}\in\mathcal{V}\label{eq:state_eq_compliant}
\end{equation}
where 
\begin{equation}
a_{s}(\mathbf{u},\mathbf{v})=\sum_{r=1,2}\int_{\Gamma_{S_{r}}}\frac{\left(\mathbf{k}_{r}\cdot\mathbf{u}\right)\left(\mathbf{k}_{r}\cdot\mathbf{v}\right)}{\left\Vert \mathbf{k}_{r}\right\Vert }\mathrm{d}s
\end{equation}
The cost functional for the mechanism design problem is defined as
\begin{equation}
J(\rho)=-\int_{\Gamma_{S_{1}}}\mathbf{k}_{1}\cdot\mathbf{u}_{\rho}\mathrm{d}s+\lambda\int_{\Omega}\rho\mathrm{d}\mathbf{x}\label{eq:Compliant-obj}
\end{equation}
where the second term again represents a constraint on the volume
of the design. The first term of this objective is a measure of the
(negative of) force applied to the workpiece in the direction of $\mathbf{k}_{1}$
which can be seen from the following relation:
\begin{equation}
\int_{\Gamma_{S_{1}}}\mathbf{k}_{1}\cdot\mathbf{u}_{\rho}\mathrm{d}s=\int_{\Gamma_{S_{1}}}\left[-\mathbf{t}_{S_{1}}(\mathbf{u}_{\rho})\right]\cdot\frac{\mathbf{k}_{1}}{\left\Vert \mathbf{k}_{1}\right\Vert }\mathrm{d}s
\end{equation}
Viewed another way, the minimization of the first term of (\ref{eq:Compliant-obj})
amounts to maximizing the displacement of the structure at the location
of the workpiece in the direction of $\mathbf{k}_{1}$.{\footnotesize{}}
\begin{table}
\begin{centering}
{\footnotesize{}}%
\begin{tabular}{|c|c|c||c|c|c|c|c|c|c|c|}
\hline 
{\footnotesize{algorithm}} & {\footnotesize{$\tau_{0}$}} & {\footnotesize{$m_{n}$}} & {\footnotesize{\# it.}} & {\footnotesize{\# bt.}} & {\footnotesize{$\ell(\mathbf{u}_{\rho})$}} & {\footnotesize{$R(\rho)$}} & {\footnotesize{$V(\rho)$}} & {\footnotesize{$\tilde{J}(\rho)$}} & {\footnotesize{$E_{1}$}} & {\footnotesize{$E_{2}$}}\tabularnewline
\hline 
\hline 
{\footnotesize{TMP}} & {\footnotesize{1}} & {\footnotesize{1}} & {\footnotesize{138}} & {\footnotesize{0}} & {\footnotesize{102.306}} & {\footnotesize{4.669}} & {\footnotesize{0.4740}} & {\footnotesize{201.779}} & {\footnotesize{6.989e-6}} & {\footnotesize{1.978e-5}}\tabularnewline
\hline 
{\footnotesize{TMP}} & {\footnotesize{2}} & {\footnotesize{1}} & {\footnotesize{169}} & {\footnotesize{62}} & {\footnotesize{102.716}} & {\footnotesize{4.075}} & {\footnotesize{0.4720}} & {\footnotesize{201.189}} & {\footnotesize{9.780e-6}} & {\footnotesize{1.679e-5}}\tabularnewline
\hline 
{\footnotesize{TMP}} & {\footnotesize{1}} & {\footnotesize{0.03}} & {\footnotesize{153}} & {\footnotesize{0}} & {\footnotesize{100.738}} & {\footnotesize{5.185}} & {\footnotesize{0.4855}} & {\footnotesize{203.014}} & {\footnotesize{7.217e-6}} & {\footnotesize{1.998e-4}}\tabularnewline
\hline 
{\footnotesize{TMP}} & {\footnotesize{2}} & {\footnotesize{0.03}} & {\footnotesize{98}} & {\footnotesize{0}} & {\footnotesize{100.568}} & {\footnotesize{5.173}} & {\footnotesize{0.4862}} & {\footnotesize{202.970}} & {\footnotesize{9.795e-6}} & {\footnotesize{1.566e-4}}\tabularnewline
\hline 
\end{tabular}
\par\end{centering}{\footnotesize \par}

{\footnotesize{\caption{{\footnotesize{Summary of the results for the MBB problem with $\beta=0.01$
\label{tab:Summary-Mbb-01}}}}
}}
\end{table}
{\footnotesize \par}

The cost functional, in the discrete setting, is given by 
\begin{equation}
J(\mathbf{z})=-\mathbf{L}^{T}\mathbf{U}(\mathbf{z})+\lambda\mathbf{z}^{T}\mathbf{v}\label{eq:Inverter-objective}
\end{equation}
where $\left[\mathbf{L}\right]_{i}=\int_{\Gamma_{S_{1}}}\mathbf{k}_{1}\cdot\mathbf{N}_{i}\mathrm{d}s$
and $\mathbf{U}(\mathbf{z})$ solves and, as before, $\mathbf{U}(\mathbf{z})$
is the solution to\textbf{ $\left[\mathbf{K}(\mathbf{z})+\mathbf{K}_{s}\right]\mathbf{U}=\mathbf{F}$}.
Here \textbf{$\mathbf{K}_{s}$} is the stiffness matrix associated
with bilinear form $a_{s}(\cdot,\cdot)$ and is independent of the
design. The gradient of $J$ can be readily computed as $\nabla J(\mathbf{z})=-\mathbf{P}^{T}\overline{\mathbf{E}}(\mathbf{z})+\lambda\mathbf{v}$
where
\begin{equation}
\left[\overline{\mathbf{E}}(\mathbf{z})\right]_{e}=p\left[\mathbf{P}\mathbf{z}\right]_{e}^{p-1}\overline{\mathbf{U}}(\mathbf{z})^{T}\mathbf{k}_{e}\mathbf{U}(\mathbf{z})\label{eq:E-interverter}
\end{equation}
and $\overline{\mathbf{U}}(\mathbf{z})$ is the solution to the \emph{adjoint}
problem
\begin{equation}
\left[\mathbf{K}(\mathbf{z})+\mathbf{K}_{s}\right]\overline{\mathbf{U}}=\mathbf{L}\label{eq:Adjoint-system}
\end{equation}
For more details on the formulation of the compliant mechanism design,
we refer the reader to \citep{Sigmund:1997p3178,Bendsoe-book}. It
is evident that $\nabla J$ can take both positive and negative values.
The main implication of this for the proposed algorithm is that the
reciprocal approximation of the cost functional is not convex and
so we cannot use its Hessian directly in the proximal term of the
quadratic model. A simple alternative that we tested is to use (\ref{eq:Hn_z_rec})
with the diagonal entries modified as
\begin{equation}
h_{k}(\mathbf{y})=\left|\frac{2}{\left[\mathbf{y}\right]_{k}}\left[\mathbf{P}^{T}\overline{\mathbf{E}}(\mathbf{y})\right]_{k}\right|\label{eq:h_k_inverster}
\end{equation}
Such an approximation has been previously explored by \citep{Groenwold:2009p2639,Groenwold:2010p2681}
and is similar in spirit to approximations in Svanberg's Method of
Moving Asymptotes \citep{Svanberg:1987p1955}. We defer a more detailed
study of suitable approximation of the Hessian for general problems
to our future work which, as illustrated in this paper, must be based
on \emph{a priori} knowledge of the cost functional.

The compliant mechanism design is known to be more prone to getting
trapped in suboptimal local minima. One such local minimum is $\rho_{h}\equiv\delta_{\rho}$
where the entire structure is eliminated and (virtually) no work is
transferred between the input actuator and the output location. For
this case, the value of the cost functional is roughly zero since
there is no density variation and minimum volume of material. To avoid
converging to this solution, we use a smaller step size parameter
$\tau_{n}=0.1$. Also we begin with small volume penalty parameter
of $\lambda=0.02$ which is then increased to $\lambda=0.15$ once
the value of cost functional reaches a negative value. This point
roughly corresponds to an intermediate density distribution in which
the structure connects the input force to the output location. The
final solution for $\beta=3\times10^{-4}$, a grid of $160\times160$
quadrilateral elements, and the two-metric projection algorithm is
shown in Figure \ref{fig:Mechanism}. This solution required a total
of 140 iterations.

\section{Discussion and Concluding Remarks}

Since the splitting algorithm presented here is a first-order method,
it is also appropriate to compare its performance to the gradient
projection algorithm, which is among the most basic first-order methods
for solving constrained optimization problems. The next iterate in
the gradient projection method is simply the projection of the unconstrained
gradient descent step onto the admissible space. In the absence of
move limits and in the discrete setting, we have the following update
expression
\begin{equation}
\mathbf{z}_{n+1}=\underset{\mathbf{z}\in\left[\delta_{\rho},1\right]^{m}}{\mbox{argmin}}\ \left\Vert \mathbf{z}-\left[\mathbf{z}_{n}-\frac{\tau_{n}}{\alpha}\nabla\tilde{J}(\mathbf{z}_{n})\right]\right\Vert ^{2}\label{eq:Grad_proj_z_n+1}
\end{equation}
where the scaling parameter $\alpha=4\lambda A$ is defined as before
in order to allow for a direct comparison with the forward-backward
splitting in the case $\mathbf{H}_{n}=\alpha\mathbf{I}$. We determine
the step size parameter $\tau_{n}$ in each iteration using the backtracking
procedure (\ref{eq:Backtracking}) based on the Armijo-type descent
condition (\ref{eq:Armijo}). Note that due to the simple structure
of the constraint set, computing the gradient $\nabla\tilde{J}$ constitutes
the main computational cost of the gradient projection algorithm in
each iteration. Table \ref{tab:Summary-gradproj-MMA} summarizes the
results for the MBB beam problem with $\beta=0.06$ for two different
choice of initial step size parameter $\tau_{0}$. First observe that
the step sizes are smaller compared to the forward-backward algorithm,
a fact that can be seen from the equivalent expression for (\ref{eq:Grad_proj_z_n+1})
given by
\begin{equation}
\mathbf{z}_{n+1}=\underset{\mathbf{z}\in\left[\delta_{\rho},1\right]^{m}}{\mbox{argmin}}\ \tilde{J}(\mathbf{z}_{n})+\left(\mathbf{z}-\mathbf{z}_{n}\right)^{T}\nabla\tilde{J}(\mathbf{z}_{n})+\frac{1}{2\tau_{n}}\left\Vert \mathbf{z}-\mathbf{z}_{n}\right\Vert _{\alpha\mathbf{I}}^{2}
\end{equation}
This shows that in each iteration, we construct a quadratic model
for the composite objective $\tilde{J}$. By constrast, the quadratic
model in (\ref{eq:z-Quad-Subproblem}) is only used for $J$ and the
regularization term appears exactly. Since $\nabla\tilde{J}$ has
a larger Lipschitz constant compared to $\nabla J$, it is therefore
expected that $\tau_{n}$ must be smaller to ensure descent. It is
also instructive to recall the informal derivation of the forward-backward
algorithm in \citep{Talischi-JMPS} where the main difference with
the gradient projection algorithm was the use of a semi-implicit (in
place of an explicit) temporal discertization of the gradient flow
equation. Note that the gradient projection algorithm converged to
the same solution as before (cf. Figure \ref{fig:Mbb-result-06})
though in the case of $\tau_{0}=0.25$, the convergence was too slow
and we terminated the algorithm after 1,000 iterations. 

Since the Method of Moving Asymptotes \citep{Svanberg:1987p1955}
is the most widely used algorithm in the topology optimization literature,
we also tested its performance using the same MBB problem. We followed
the common practice and used the algorithm as a black-box optimization
routine. In particular, we provided the algorithm with the gradient
of composite objective $\tilde{J}$ and did not make any changes to
the open source code provided by Svangberg%
\footnote{We remark that with a few exceptions, MMA is used in the same way
by Borrvall in a review paper \citep{Borrvall:2001p2519} where he
compares various regularizations schemes, including Tikhonov regularization.%
}. MMA internally generates a separable convex approximation to $\tilde{J}$
using reciprocal-type expansions with appropriately defined and updated
asymptotes. Though such approximations are suitable for the structural
term, they may be inaccurate for the Tikhonov regularizer and the
composite objective. As shown in Table \ref{tab:Summary-gradproj-MMA},
MMA did not converge (according to the convergence criteria described
earlier) in 1,000 iterations before it was terminated. Furthermore,
not only was the final value of the objective function larger than
that obtained by gradient projection or either splitting algorithm,
the final density was topologically different from the solution shown
in the Figure \ref{fig:Mbb-result-06}.

The fact that the present splitting framework outperforms MMA should
not be surprising. Unlike MMA, which is far more general and can handle
a much broader class of problems \citep{Svanberg:2001p1960}, the
present algorithm is tailored to the specific structure of (\ref{eq:reg-compliance})
(or (\ref{eq:Discrete-opt}) in the discrete setting) and provides
an ideal treatment of its constituents. First, the composite objective
is the sum of two terms and algorithm deals with each term separately.
The regularization term $R$ is represented with a high degree fidelity
since the resulting subproblem with its simple structure can be solved
efficiently. The structural term $J$, while expensive to compute,
contains many local minima and very fast convergence usually at best
reaches a suboptimal local minimum. Moreover, $J$ tends to be rather
flat near stationary points and so one should not require a high level
accuracy for satisfaction of the first order conditions of optimality.
As a side remark, these characteristics indicate that second order
methods do not pay off given their significantly higher computational
cost per iteration%
\footnote{Computing the exact Hessian information is especially expensive for
PDE-constrained problem since every Hessian-vector product requires
the solution to an adjoint system.%
}. The other drawback of using exact second order information is the
storage requirements, quadratic in the size of the problem, which
can be prohibitive for large-scale problems such as those encountered
in practical applications of topology optimization. Therefore, first
order methods are better suited for minimizing $J$. 

In the splitting algorithm proposed here, we use additional knowledge
about the behavior of $J$ to construct accurate approximations using
only first order information and minimal storage requirements. Furthermore,
the two-metric approach allows for a computationally efficient treatment
of the constraint set. In fact, the proposed approach is aligned with
the renewed interest in first-order convex optimization algorithms
for solving large-scale inverse problems in signal recovery, statistical
estimation, and machine learning \citep{SJWrigthMachine,Combettes:2006p3658,Bredies:2009p3671,Duchi:2009p3686}.
Our rather restricted and narrow comparison with MMA is meant to motivate
the virtue of developing such tailored algorithms. We note that, aside
from efficiency, robustness is also a major issue for solving topology
optimization problems (see, for example, comments in \citep{Borrvall:2001p2519}
on total variation regularization). Although the high sensitivity
to parameters is, to a large extent, intrinsic to the size, nonconvexity
and sometimes nonsmoothness of these problems, we emphasize that it
should be minimized as much as possible. Developing an appropriately-designed
optimization algorithm that fits the structure of the problem at hand
can be key to achieving this. 
\begin{table}
\begin{centering}
{\footnotesize{}}%
\begin{tabular}{|c|c|c|c|c|c|}
\hline 
{\footnotesize{algorithm}} & {\footnotesize{$\tau_{0}$}} & {\footnotesize{\# it.}} & {\footnotesize{\# bt.}} & {\footnotesize{$\tilde{J}(\rho)$}} & {\footnotesize{$E_{2}$}}\tabularnewline
\hline 
\hline 
{\footnotesize{GP}} & {\footnotesize{0.25}} & {\footnotesize{1000{*}}} & {\footnotesize{0}} & {\footnotesize{210.74}} & {\footnotesize{1.362e-4}}\tabularnewline
\hline 
{\footnotesize{GP}} & {\footnotesize{0.5}} & {\footnotesize{568}} & {\footnotesize{79}} & {\footnotesize{210.68}} & {\footnotesize{8.939e-5}}\tabularnewline
\hline 
{\footnotesize{MMA}} & -- & {\footnotesize{1000{*}}} & {\footnotesize{0}} & {\footnotesize{213.39}} & {\footnotesize{1.913e-4}}\tabularnewline
\hline 
\end{tabular}
\par\end{centering}{\footnotesize \par}

\caption{Summary of the results for gradient projection and MMA algorithm for
the MBB beam problem with $\beta=0.06$. The asterisk indicates that
the maximum allowed iteration count of 1,000 was reached before the
convergence criteria was met \label{tab:Summary-gradproj-MMA}}
\end{table}

In the extensions of this work, we intend to consider nonsmooth regularizers
such as the total variation of density within the present variable
metric scheme. This would require the extension of available denoising
algorithms (e.g. \citep{Chambolle:2004p2741,Bredies:2009p3671}) for
solving the resulting subproblems in each iteration. Also of interest
is the use of accelerated first order methods such as those proposed
in \citep{Nesterov2007} and \citep{Beck2008} that can improve the
convergence speed of the algorithms. Developing a two-metric variation
of such algorithms for the constrained minimization problems of topology
optimization is promising.

\section*{Acknowledgements}

The authors acknowledge the support by the Department of Energy Computational
Science Graduate Fellowship Program of the Office of Science and National
Nuclear Security Administration in the Department of Energy under
contract DE-FG02-97ER25308.

\bibliographystyle{siam}
\bibliography{bib-May2012-thesis}

\begin{thebibliography}{10}

\bibitem{Adlers-thesis}
{\sc M.~Adlers}, {\em Sparse Least Squares Problems with Box Constraints},
  Department of Mathematics, Linkoping University, Thesis, 1998.

\bibitem{Allaire-book}
{\sc G.~Allaire}, {\em Shape Optimization by the Homogenization Method},
  Springer, New York, 2001.

\bibitem{Arora:1980p3660}
{\sc J.~S. Arora}, {\em Analysis of optimality criteria and gradient projection
  methods for optimal structural design}, Comput Methods Appl Mech Engrg, 23
  (1980), pp.~185--213.

\bibitem{Beck2008}
{\sc A.~Beck and M.~Teboulle}, {\em A fast iterative shrinkage-thresholding
  algorithm for linear inverse problems}, SIAM J Image Sci, 2 (2008),
  pp.~183--202.

\bibitem{Beck-review}
\leavevmode\vrule height 2pt depth -1.6pt width 23pt, {\em Gradient-based
  algorithms with applications to signal recovery problems}, in Convex
  Optimization in Signal Processing and Communications, Cambridge university
  press, 2010.

\bibitem{Bendsoe:1989p447}
{\sc M.~P. Bends{\o}e}, {\em Optimal design as material distribution probelm},
  Struct Optimization, 1 (1989), pp.~193--202.

\bibitem{Bendsoe-book}
{\sc M.~P. Bends{\o}e and O.~Sigmund}, {\em Topology Optimization: Theory,
  {M}ethods and {A}pplications}, Springer, 2003.

\bibitem{Bert-ProjNewton}
{\sc D.~P. Bertsekas}, {\em Projected newton methods for optimization problems
  with simple constraints}, SIAM J Control Opt, 20 (1982), pp.~221--246.

\bibitem{Bert-book}
\leavevmode\vrule height 2pt depth -1.6pt width 23pt, {\em Nonlinear
  Programming}, Athena Scientific, 2nd~ed., 1999.

\bibitem{Borrvall:2001p2519}
{\sc T.~Borrvall}, {\em Topology optimization of elastic continua using
  restriction}, Arch Comput Method E, 8 (2001), pp.~251--285.

\bibitem{Borrvall:2001p422}
{\sc T.~Borrvall and J.~Petersson}, {\em Topology optimization using
  regularized intermediate density control}, Comput Methods Appl Mech Engrg,
  190 (2001), pp.~4911--4928.

\bibitem{Bourdin:2001p457}
{\sc B.~Bourdin}, {\em Filters in topology optimization}, Int J Numer Meth Eng,
  50 (2001), pp.~2143--2158.

\bibitem{Bourdin:2003p1799}
{\sc B.~Bourdin and A.~Chambolle}, {\em Design-dependent loads in topology
  optimization}, {ESAIM} Contr Optim Ca, 9 (2003), pp.~19--48.

\bibitem{Bredies:2009p3671}
{\sc K.~Bredies}, {\em A forward--backward splitting algorithm for the
  minimization of non-smooth convex functionals in {B}anach space}, Inverse
  Probl, 25 (2009), p.~015005.

\bibitem{Brenner-FE}
{\sc S.~C. Brenner and L.~R. Scott}, {\em The Mathematical Theory of Finite
  Element Methods}, Springer, 2nd~ed., 2002.

\bibitem{Bruns:2001p2341}
{\sc T.~Bruns and D.~A. Tortorelli}, {\em Topology optimization of non-linear
  elastic structures and compliant mechanisms}, Comput Methods Appl Mech Engrg,
  190 (2001), pp.~3443--3459.

\bibitem{Burger:2006p1864}
{\sc M.~Burger and R.~Stainko}, {\em Phase-field relaxation of topology
  optimization with local stress constraints}, SIAM J Control Optim, 45 (2006),
  pp.~1447--1466.

\bibitem{Chambolle:2004p2741}
{\sc A.~Chambolle}, {\em An algorithm for total variation minimization and
  applications}, J Math Imaging Vis, 20 (2004), pp.~89--97.

\bibitem{Chen:1997p3679}
{\sc G.~H.~G. Chen and R.~T. Rockafellar}, {\em Convergence rates in
  forward-backward splitting}, SIAM J Optimiz, 7 (1997), pp.~421--444.

\bibitem{Cohen:1978p3656}
{\sc G.~Cohen}, {\em Optimization by decomposition and coordination: a unified
  approach}, IEEE Trans Autom Control, 23 (1978), pp.~222--232.

\bibitem{Combettes:2006p3658}
{\sc P.~L. Combettes and V.~R. Wajs}, {\em Signal recovery by proximal
  forward-backward splitting}, Multiscale Model Sim, 4 (2006), pp.~1168--1200.

\bibitem{Duchi:2009p3686}
{\sc J.~Duchi and Y.~Singer}, {\em Efficient online and batch learning using
  forward backward splitting}, J Mach Learn Res, 10 (2009), pp.~2899--2934.

\bibitem{Evans}
{\sc L.~C. Evans}, {\em Partial Differential Equations}, Graduate Studies in
  Mathematics, American Mathematical Society, Rhode Island, 1998.

\bibitem{Gafni-TMP}
{\sc E.~M. Gafni and D.~Bertsekas}, {\em Two-metric projection methods for
  constrained optimization}, SIAM J Control Opt, 20 (1984), pp.~936--964.

\bibitem{Groenwold:2008p2685}
{\sc A.~A. Groenwold and L.~F.~P. Etman}, {\em On the equivalence of optimality
  criterion and sequential approximate optimization methods in the classical
  topology layout problem}, Int J Numer Meth Eng, 73 (2008), pp.~297--316.

\bibitem{Groenwold:2009p2639}
\leavevmode\vrule height 2pt depth -1.6pt width 23pt, {\em A quadratic
  approximation for structural topology optimization}, Int J Numer Meth Eng, 82
  (2010), pp.~505--524.

\bibitem{Groenwold:2010p2681}
{\sc A.~A. Groenwold, L.~F.~P. Etman, and D.~W. Wood}, {\em Approximated
  approximations for {SAO}}, Struct Multidisc Optim, 41 (2010), pp.~39--56.

\bibitem{Nesterov2007}
{\sc Y.~Nesterov}, {\em Gradient methods for minimizing composite objective
  function}.
\newblock available at http://www.ecore.be/DPs/dp1191313936.pdf., 2007.

\bibitem{Patriksson:1998p3659}
{\sc M.~Patriksson}, {\em Cost approximation: A unified framework of descent
  algorithms for nonlinear programs}, SIAM J Optimiz, 8 (1998), pp.~561--582.

\bibitem{Petersson:1999p452}
{\sc J.~Petersson}, {\em Some convergence results in perimeter-controlled
  topology optimization}, Comput Methods Appl Mech Engrg, 171 (1999),
  pp.~123--140.

\bibitem{Petersson:1998p440}
{\sc J.~Petersson and O.~Sigmund}, {\em Slope constrained topology
  optimization}, Int. J. Numer. Meth. Engng, 41 (1998), pp.~1417--1434.

\bibitem{Rozvany:2009p2274}
{\sc G.~I.~N. Rozvany}, {\em A critical review of established methods of
  structural topology optimization}, Struct Multidisc Optim, 37 (2009),
  pp.~217--237.

\bibitem{Rozvany:1992p981}
{\sc G.~I.~N. Rozvany, M.~Zhou, and T.~Birker}, {\em Generalized shape
  optimization without homogenization}, Struct Optimization, 4 (1992),
  pp.~250--252.

\bibitem{Sigmund:1997p3178}
{\sc O.~Sigmund}, {\em On the design of compliant mechanisms using topology
  optimization}, Mechanics Based Design of Structures and Machines, 25 (1997),
  pp.~493--524.

\bibitem{Sigmund:2007p460}
\leavevmode\vrule height 2pt depth -1.6pt width 23pt, {\em Morphology-based
  black and white filters for topology optimization}, Struct Multidisc Optim,
  33 (2007), pp.~401--424.

\bibitem{SigmundMaute2013}
{\sc O.~Sigmund and K.~Maute}, {\em Sensitivity filtering from a continuum
  mechanics perspective}, Struct Multidisc Optim, 46 (2012), pp.~471--475.

\bibitem{Sigmund:1998p441}
{\sc O.~Sigmund and J.~Petersson}, {\em Numerical instabilities in topology
  optimization: A survey on procedures dealing with checkerboards,
  mesh-dependencies and local minima}, Struct Optimization, 16 (1998),
  pp.~68--75.

\bibitem{Svanberg:1987p1955}
{\sc K.~Svanberg}, {\em The {M}ethod of {M}oving {A}symptotes--{A} new method
  for structural optimization}, Int J Numer Meth Eng, 24 (1987), pp.~359--373.

\bibitem{Svanberg:2001p1960}
\leavevmode\vrule height 2pt depth -1.6pt width 23pt, {\em A class of globally
  convergent optimization methods based on conservative convex separable
  approximations}, SIAM J Optimiz, 12 (2001), pp.~555--573.

\bibitem{Takezawa:2010p2730}
{\sc A.~Takezawa, S.~Nishiwaki, and M.~Kitamura}, {\em Shape and topology
  optimization based on the phase field method and sensitivity analysis}, J
  Comput Phys, 229 (2010), pp.~2697--2718.

\bibitem{Talischi-JMPS}
{\sc C.~Talischi and G.~H. Paulino}, {\em An operator splitting algorithm for
  {T}ikhonov-regularized topology optimization}, Comput Methods Appl Mech
  Engrg, 253 (2013), pp.~599--608.

\bibitem{Talischi:2010p3179}
{\sc C.~Talischi, G.~H. Paulino, A.~Pereira, and I.~F.~M. Menezes}, {\em
  Polygonal finite elements for topology optimization: A unifying paradigm},
  Int J Numer Meth Eng, 82 (2010), pp.~671--698.

\bibitem{PolyTop}
\leavevmode\vrule height 2pt depth -1.6pt width 23pt, {\em {PolyTop}: a
  {M}atlab implementation of a general topology optimization framework using
  unstructured polygonal finite element meshes}, Struct Multidisc Optim, 45
  (2012), pp.~329--357.

\bibitem{SJWrigthMachine}
{\sc S.~J. Wright}, {\em Optimization in machine learning}, in Neural
  Information Processing Systems (NIPS) Workshop, 2008.

\end{thebibliography}

\end{document}